\documentclass[12pt]{amsart}

\usepackage{amsmath, amsthm, amsfonts, amssymb}
\usepackage{epsfig, enumerate}
\usepackage{color}

\usepackage{hyperref}
\usepackage[all]{xy}

\usepackage{float}
\usepackage[left=3.5cm, marginratio=1:1]{geometry}

\usepackage{tikz-cd}
\usepackage{tikz}
\usetikzlibrary{arrows,shapes,trees}
\usepackage{graphicx}
\usepackage{amssymb}
\usepackage{amsfonts}
\usepackage{amsthm}
\usepackage{amsmath}
\usepackage{mathrsfs}
\usepackage{mathtools}
\usepackage{extpfeil}

\newtheorem{maintheorem}{Theorem} 

\newtheorem{maintheoremprime}{Theorem}

\newtheorem{theorem}{Theorem}[section]
\newtheorem{lemma}[theorem]{Lemma}
\newtheorem{corollary}[theorem]{Corollary}
\newtheorem{proposition}[theorem]{Proposition}

\newtheorem{remark}[theorem]{Remark}
\newtheorem{example}[theorem]{Example}

\numberwithin{equation}{section}

\newcommand {\ignore}[1]  {}

\newcommand{\CC}{{\mathbb C}}

\newcommand{\NN}{{\mathbb N}}

\newcommand{\RR}{{\mathbb R}}

\newcommand{\TT}{{\mathbb T}}

\newcommand{\ZZ}{{\mathbb Z}}

\newcommand{\cA}{{\mathcal A}}

\newcommand{\cC}{{\mathcal C}}

\newcommand{\cE}{{\mathcal E}}

\newcommand{\cG}{{\mathcal G}}
\newcommand{\cH}{{\mathcal H}}

\newcommand{\cL}{{\mathcal L}}

\newcommand{\cO}{{\mathcal O}}
\newcommand{\cP}{{\mathcal P}}

\newcommand{\cS}{{\mathcal S}}

\newcommand{\cU}{{\mathcal U}}

\newcommand{\cZ}{{\mathcal Z}}

\newcommand{\fT}{{\mathfrak T}}

\newcommand{\fX}{{\mathfrak X}}
\newcommand{\fY}{{\mathfrak Y}}

\newcommand{\fg}{{\mathfrak g}}

\newcommand{\fk}{{\mathfrak k}}

\newcommand{\fn}{{\mathfrak n}}

\newcommand{\fr}{{\mathfrak r}}

\newcommand{\Ad}{{\rm Ad}}
\newcommand{\ad}{{\rm ad}}
\newcommand{\id}{{\rm id}}
\newcommand{\Lip}{{\rm Lip}}
\newcommand{\reg}{{\rm reg}}
\newcommand{\adj}{{\rm adj}}
\newcommand{\Aut}{{\rm Aut}}
\newcommand{\diag}{{\rm diag}}
\newcommand{\Hom}{{\rm Hom}}
\newcommand{\End}{{\rm End}}
\newcommand{\Aff}{{\rm Aff}}

\newcommand{\SL}{{\rm SL}}
\newcommand{\SO}{{\rm SO}}
\newcommand{\Lie}{{\rm Lie}}
\newcommand{\GL}{{\rm GL}}
\newcommand{\tr}{{\rm tr}}

\newcommand{\vol}{{\rm vol}}
\newcommand{\Sym}{{\rm Sym}}
\newcommand{\Diff}{{\rm Diff}}
\newcommand{\dist}{{\rm dist}}
\newcommand{\Haar}{{\rm Haar}}
\newcommand{\Isom}{{\rm Isom}}
\newcommand{\supp}{{\rm supp}}
\newcommand{\Span}{{\rm Span}}

\begin{document}

\title[Ergodicity and Symmetry]{On ergodicity and symmetry for affine systems}

\author{Danijela Damjanovi\'c}

\address[Damjanovi\'c]{Department of mathematics, Kungliga Tekniska högskolan, Lindstedtsv\"agen 25, SE-100 44 Stockholm, Sweden.} 
\email{ddam@kth.se}

\author{Amie Wilkinson}
\address[Wilkinson]{Department of mathematics, the University of Chicago, Chicago, IL, US, 60637}
\email{wilkinso@math.uchicago.edu}

\author{Chengyang Wu}
\address[Wu]{Department of Mathematics, University of Chicago, Chicago, IL, US, 60637}
\email{chengyangwu1999@gmail.com}

\author{Disheng Xu}
\address[Xu]{School of Sciences, Great Bay University and Great Bay Institute for Advanced Study, Songshan Lake Innovation Entrepreneurship Community A5,  Dongguan, Guangdong, China, 523000}
\email{xudisheng@gbu.edu.cn}

\begin{abstract}
Let $f$ be a volume-preserving affine diffeomorphism of a compact homogeneous space with smooth centralizer $\cZ^{\infty}(f)$. We establish the following two equivalences: \begin{equation*}
\begin{aligned}
    f\text{ is ergodic}&\Longleftrightarrow \cZ^{\infty}(f)\text{ is a finite-dimensional Lie group};\\
    f\text{ is non-ergodic}&\Longleftrightarrow \cZ^{\infty}(f)\text{ contains either }C_c^{\infty}((0,1),\RR) \text{ or }\Diff_c^{\infty}((0,1)).
\end{aligned}
\end{equation*}
This strengthens the corresponding perturbative version established in \cite{DWWX}.

We also prove analogues of the above equivalences for affine-induced diffeomorphisms on bi-homogeneous spaces, frame flows over negatively curved manifolds, and automorphisms of Kummer-type complex surfaces.
\end{abstract}

\maketitle

\tableofcontents

\section{Introduction} 
This paper focuses on the following fundamental question: {\em When is the group of symmetries of a diffeomorphism a Lie group?}
Analogous questions arise in other contexts: for example,  
all isometries of a complete Riemannian manifold form a Lie group \cite{MS39}, and the biholomorphic automorphisms of the Riemann sphere are the M\"obius group $\mathrm{PSL}(2,\mathbb C)$. 


In the dynamical setting, the symmetry group of a diffeomorphism $f$ of a closed manifold $M$ is its {\em smooth centralizer:}
\[
    \cZ^\infty(f)
    :=\{h\in\Diff^\infty(M):h\circ f=f\circ h\},
\]
equipped with the topology induced by $\Diff^\infty(M)$.
This group always contains the iterates of $f$, but its additional
symmetries can vary considerably. At one extreme,
$\cZ^\infty(\id_M)=\Diff^\infty(M)$; at the other extreme,
the centralizer of a diffeomorphism may be {\em trivial}, in the sense that it consists only of its iterates.
We study here which dynamical properties force the centralizer to be a
finite-dimensional Lie group.

The study of centralizers has a long history, including Kopell's work \cite{Kopell} on commuting diffeomorphisms of the circle and the rigidity results of Palis--Yoccoz \cite{PY89}. One central theme is Smale's question concerning generic triviality of centralizers, answered in the $C^1$ topology by Bonatti, Crovisier and Wilkinson \cite{BCW}. A complementary problem is to understand systems with additional symmetries and to identify the dynamical mechanisms that constrain the structure of their centralizers.

In our previous work \cite{DWWX}, we studied centralizers of affine diffeomorphisms on compact homogeneous spaces. We proved that ergodicity implies that the smooth centralizer is a Lie group, and established a perturbative converse. Our first main result, Theorem~\ref{T:main}, completes the characterization for individual affine diffeomorphisms:
\[  f\text{ is ergodic}
    \quad\Longleftrightarrow\quad
    \cZ^\infty(f)\text{ is a finite-dimensional Lie group}.
\]
Here ergodicity is taken with respect to the 
Haar probability measure. In the non-ergodic case, we exhibit
infinite-dimensional subgroups of the centralizer of $f$ itself.

The relation between ergodicity and rigidity of symmetries has substantial precedents in homogeneous dynamics.
Parry \cite{Parry71} proved that measure-preserving conjugacies between ergodic unipotent affine transformations of compact nilmanifolds are affine almost everywhere.
Walters \cite{Wa-example,Wa} showed that the topological centralizer of an ergodic nilmanifold automorphism is affine and obtained the same conclusion for weakly mixing affine transformations.
Witte \cite{Witte87,Wit87cor} extended this measurable affine rigidity to ergodic unipotent affine transformations and also established it for weakly mixing affine transformations with zero entropy.

In these results,  weak mixing implies that symmetries must preserve the given affine structure on $G/B$, forming not just a Lie group, but a subgroup of the affine group $\Aff(G/B)$. Ergodicity alone does not imply this stronger conclusion: Walters \cite{Wa-example} constructed ergodic affine transformations of tori with non-affine smooth symmetries. The Lie group conclusion in our results therefore captures a form of rigidity that persists beyond affine rigidity.

A central ingredient in our proof is a regularity property of invariant functions: for an affine diffeomorphism, we show that smooth (even analytic) invariant functions are dense in the space of square-integrable invariant functions (Theorem~\ref{T:aux-main}). In particular,  non-ergodicity supplies a nonconstant smooth invariant function. We develop constructions that combine such functions with invariant vector fields to produce infinite-dimensional families of commuting diffeomorphisms.

Our arguments build on the classical ergodic theory of homogeneous flows. Our proof uses Moore-Mautner's theorem \cite{Moore80} and a modification/correction of Dani's suspension construction \cite{Dani77}. The existence of smooth invariant functions for non-ergodic systems was established by Brezin--Moore in the context of homogeneous flows \cite[Theorem 6.1]{BM81}.
Theorem~\ref{T:aux-main} establishes the stronger {\em density} statement for affine diffeomorphisms. Related geometric descriptions of the ergodic decomposition of homogeneous flows were obtained by Starkov \cite{Starkov91}.

These methods also apply beyond homogeneous spaces. For example, we prove analogous results for affine-induced diffeomorphisms of bihomogeneous spaces (Theorem \ref{T:bi-homo}) and frame flows in negative curvature (Theorem \ref{T:frame-flow}).
(These settings overlap in the locally symmetric case: frame flows on compact quotients of negatively curved symmetric spaces admit bihomogeneous models.) In both cases, we again obtain an equivalence between ergodicity and Lie centralizer. 
A related geometric result is due to Mitsutani \cite{Mitsutani}, who proved that the foliated centralizer of the geodesic flow on the universal cover of a closed negatively curved manifold is a Lie group when its stable and unstable foliations are $C^1$.

We also establish an equivalence between ergodicity and Lie centralizer for automorphisms of Kummer-type complex surfaces (Theorem~\ref{T:c-surface}).
In complex dynamics, strong restrictions on centralizers are already known: Blanc--Cantat \cite{BlancCantat} proved that the cyclic group generated by a positive-entropy automorphism of a complex projective surface has finite index in its birational centralizer. Our results concern the full smooth centralizer, without requiring its elements to preserve the complex structure.


\subsection{Main results on affine ergodic systems and Lie centralizers}
We first present our main theorem on the equivalence between ergodicity and Lie centralizer for affine systems on compact homogeneous spaces.

Let $X=G/B$ be a compact homogeneous manifold, where $G$ is a connected Lie group acting transitively and locally effectively \footnote{An action of a group $G$ on a manifold $M$ is \textsl{effective} (resp. \textsl{locally effective}) if the elements that fix every point in $M$ form a trivial (resp. discrete) subgroup.} on $X$, and $B$ is the point stabilizer of this action. We always assume that $X$ admits a $G$-invariant probability measure $\mu_{\Haar}$.

For the compact homogeneous manifold $X$, we choose and fix an \textsl{effective presentation} $(G,B)$, that is, we write $X=G/B$ where $B$ contains no nontrivial normal subgroup of $G$. With respect to this presentation, an \textsl{aﬀine diffeomorphism} is a map $f\colon X\to X$ of the form $L_a\circ \Phi$, where $L_a$ denotes left translation by $a\in G$ and $\Phi\colon G\to G$ is an automorphism preserving $B$. We denote by $\Aff(G,B)$ the space of all such aﬀine diffeomorphisms, endowed with the $C^0$ topology. When there is no ambiguity, we also write $\Aff(X)=\Aff(G,B)$ for the sake of simplicity.

From the perspective of dynamics, any $f\in\Aff(X)$ defines a measure-preserving system $(X,\mu_{\Haar},f)$. We consider the ergodicity of such system and put $$
\cE(X):=\{f\in \Aff(X):f\text{ is ergodic (with respect to }\mu_{\Haar})\}.
$$

From the perspective of geometry, 
we investigate the centralizer of $f\in\Aff(X)$ in the full group $\Diff^{\infty}(X)$, which describes its symmetries. We focus on the set $$
\cL(X):=\{f\in \Aff(X):\cZ^{\infty}(f)\text{ is a (finite-dimensional) Lie group}\}.
$$

Our main result builds the bridge between ergodicity and symmetry for affine systems:

\begin{maintheorem}\label{T:main}
For any compact homogeneous space $X$, we have $$
\cE(X)=\cL(X).
$$
Moreover, when $f\in\Aff(X)$ is non-ergodic, $\cZ^{\infty}(f)$ contains a continuously immersed copy of one of the following groups:
\begin{itemize}
    \item the abelian group $C^{\infty}_c((0,1),\RR)$ of compactly-supported smooth functions on the interval $(0,1)$ (under addition); or
    \item the simple group $\Diff_c^{\infty}((0,1))$ of compactly-supported diffeomorphisms of the interval $(0,1)$ (under composition).
\end{itemize}
\end{maintheorem}

The inclusion $\cE(X)\subseteq\cL(X)$ is proved in \cite[Theorem B]{DWWX}. For the converse inclusion $\cE(X)\supseteq \cL(X)$, we present two proofs. The first is a streamlined version of an argument that was developed entirely by the authors. The second proof was suggested by an AI assistant with which we had previously discussed the first proof.

The first, longer proof carries more information about the dynamics of  non-ergodic affine diffeomorphisms. It is closer in spirit to our earlier, more geometric,  work on this subject, which also lends itself to the study of frame flows in a nonhomogeneous setting. 

The second proof is quite short and extends to bi-homogeneous diffeomorphisms (such as geodesic and time-$T$ maps of $k$-frame flows of locally symmetric manifolds). Moreover, it yields a stronger conclusion that the smooth centralizer always contains a subgroup corresponding to the first possibility in the non-ergodic case of Theorem \ref{T:main}, namely, one isomorphic to $C^{\infty}_c((0,1),\RR)$.

\begin{maintheoremprime}\label{T:main-prime}
Let $X$ be a compact homogeneous space and let $f$ be an affine diffeomorphism on $X$. Then $f$ is non-ergodic if and only if $\cZ^{\infty}(f)$ contains a continuously immersed copy of $C^{\infty}_c((0,1),\RR)$.
\end{maintheoremprime}

We emphasize that both proofs use the key fact below 
that the smooth invariant functions are dense in the space of all square integrable invariant functions for affine systems, 
which is of independent interest:

\begin{maintheorem}\label{T:aux-main}
Let $X$ be a compact homogeneous space and let $f$ be an affine diffeomorphism on $X$. Then the subspace $C^{\infty}(X)^f$ of all $f$-invariant smooth functions on $X$ is dense in the space $L^2(X)^f$ of all $f$-invariant $L^2$ functions on $X$.

    In particular, $f$ is non-ergodic if and only if there exists a nonconstant $f$-invariant smooth function on $X$.
\end{maintheorem}

We remark that for non-ergodic homogeneous flows, the existence of a nonconstant invariant function was proved in \cite[Theorem 6.1]{BM81} by analyzing the maximal solvable and semisimple quotients. In this paper, we give an independent proof based on the density statement in Theorem \ref{T:aux-main}, and then extend it to non-ergodic affine maps using a modification of Dani’s suspension.

 It is interesting to contrast Theorem~\ref{T:aux-main} with work of Flaminio, Forni and Rodriguez Hertz \cite{FFRH}, who proved that an affine diffeomorphism of a connected finite volume homogeneous space admits infinitely many linearly independent invariant distributions unless it is smoothly conjugate to an ergodic toral translation.


\begin{remark}
One may also extend our main Theorem (Theorem \ref{T:main}) to the setting of conjugacies between two affine systems and the joinings of such systems. More precisely, on a connected compact homogeneous space, any Lipschitz conjugacy between an ergodic affine diffeomorphism and another affine diffeomorphism must be smooth, and all such conjugacies
form a finite-dimensional family. If the first map is non-ergodic and at least one smooth conjugacy exists, there are arbitrarily large finite-dimensional families of smooth conjugacies.

To connect conjugacies with joinings, observe that a conjugacy $h\colon X\to X$ between two affine diffeomorphisms $f_1,f_2$ gives an invariant graph measure:
\[
h\circ f_1=f_2\circ h\quad\Longrightarrow\quad
(f_1\times f_2)_*(\operatorname{id},h)_*\mu_{\Haar}=(\operatorname{id},h)_*\mu_{\Haar}.
\]
This measure has marginals $\mu_{\Haar}$ and $h_*\mu_{\Haar}$. When $f_1$ is ergodic, the graph measures obtained from Lipschitz conjugacies are Haar joinings and form a finite-dimensional family. When $f_1$ is nonergodic and a smooth conjugacy exists, varying $h$ among smooth conjugacies gives arbitrarily large finite-dimensional families of invariant graph measures.
\end{remark}

\subsection{Applications to more non-homogeneous spaces}

In this subsection, we establish the equivalence between ergodicity and Lie centralizer in other settings, for example, for the affine-induced diffeomorphisms on bi-homogeneous spaces, frame flows over negatively curved manifolds, and automorphisms of Kummer-type complex surfaces. Their proofs substantially generalize the arguments in the proof of Theorem \ref{T:main}.

\subsubsection{Affine-induced maps on bi-homogeneous spaces} Let $G$ be a connected Lie group, and let $B\leq G$ be a closed cocompact subgroup such that $X=G/B$ admits a $G$-invariant probability measure $\mu_{\Haar}$. Let $K\leq G$ be a compact subgroup acting freely \footnote{An action of a group $K$ on a manifold $M$ is \textsl{free} (resp. \textsl{locally free}) if any stabilizer of this action is trivial (resp. discrete).} by left translation on $X$. This induces a bi-homogeneous space structure on $$
Y=K\backslash G/B.
$$
We write $\pi\colon X\to Y$ to be the natural projection map, and write $m=\pi_*\mu_{\Haar}$.

Let $f=L_a\circ \Phi\in\Aff(G,B)$ be such that $(c_a\circ\Phi)(K)=K$. Then it induces a diffeomorphism on $Y$ as $$
F\colon Y\to Y,\qquad KgB\mapsto Kf(g)B.
$$
As in Theorem \ref{T:main}, we build the bridge between ergodicity and symmetry for the dynamical system $(Y,m,F)$. Here we denote by $\cZ^{\Lip}(F)$ the set of bi-Lipschitz homeomorphisms of $Y$ commuting with $F$. Clearly $\cZ^{\infty}(F)\leq \cZ^{\Lip}(F)$.

\begin{theorem}\label{T:bi-homo}
    For the above system $(Y,m,F)$, we have $$
   F\text{ is ergodic with respect to }m\Longleftrightarrow \cZ^{\infty}(F) \text{ is a finite-dimensional Lie group}.
    $$
    Moreover, \begin{itemize}
        \item When $F$ is ergodic, $\cZ^{\Lip}(F)=\cZ^{\infty}(F)$ is a Lie group of dimension at most $\dim{Y}$;
        \item When $F$ is non-ergodic, $\cZ^{\infty}(F)$ contains a continuously immersed copy of $(C^{\infty}_c((0,1),\RR),+)$.
    \end{itemize}
\end{theorem}

We remark that Theorem \ref{T:bi-homo} specializes to Theorem \ref{T:main} when $K=\{1_G\}$. However, the bihomogeneous result is not a direct consequence of the homogeneous case. A reduction by lifting would first require showing that every centralizing diffeomorphism of $Y$ admits a suitable lift to $X$. Moreover, commutation downstairs only forces the commutator of the lifts to preserve each $K$-fiber; it need not be trivial or even given by a fixed element of $K$. Our proof works directly on the quotient and therefore avoids these additional lifting and compatibility problems.

\subsubsection{Frame flows over negatively curved manifolds}
Let $(N,g)$ be a connected, closed, oriented Riemannian manifold of dimension $n\geq 2$ with negative curvature. Let $FN$ denote the bundle of positively oriented orthonormal frames on $N$, and let $SN$ denote the unit tangent bundle over $N$.  
We denote by $\pi\colon FN\to SN$ the projection that keeps the first vector in each orthonormal frame. This makes $FN$ a principal $\SO(n-1)$-bundle over $SN$. Write $K=\SO(n-1)$.

For $u=(x;v_1,\ldots,v_n)\in FN$, let $\gamma_{x,v_1}$ be the geodesic with initial condition
$$
\gamma_{x,v_1}(0)=x,\qquad
\dot\gamma_{x,v_1}(0)=v_1,
$$
and let $\cP^{x,v_1}_{0,t}$ denote parallel transport along $\gamma_{x,v_1}$ from time $0$ to time $t$. The \textsl{frame flow} $\Phi=\{\Phi_t\}_{t\in\RR}$ on $FN$ is defined by
$$
\Phi_t\colon (x;v_1,\ldots,v_n)\mapsto
\left(\gamma_{x,v_1}(t);\cP^{x,v_1}_{0,t}(v_1),\cdots,\cP^{x,v_1}_{0,t}(v_n)
\right).
$$
Since parallel transport is linear, $\Phi_t$ commutes with the right action of $K$. 

Let $\mu_{\mathrm L}$ denote the normalized Liouville measure on $SN$. The natural probability measure $\mu$ on $FN$ is obtained by integrating the normalized $K$-invariant probability measures on the fibers of $\pi\colon FN\to SN$ against $\mu_{\mathrm L}$. The measure $\mu$ is smooth, has full support, and is invariant under the frame flow $\Phi$.

As in Theorem \ref{T:main}, we build the bridge between ergodicity and symmetry for the dynamical system $(FN,\mu,\Phi_T)$. Here we denote by $\cZ^0(\Phi_T)$ the set of homeomorphisms of $FN$ commuting with $\Phi_T$. Clearly $\cZ^{\infty}(\Phi_T)\subseteq\cZ^0(\Phi_T)$.

\begin{theorem}\label{T:frame-flow}
Let $T\neq 0$. For the above system $(FN,\mu,\Phi_T)$, $$
\Phi_T\text{ is ergodic with respect to }\mu\Longleftrightarrow \cZ^{\infty}(\Phi_T) \text{ is a finite-dimensional Lie group}.
$$
Moreover, \begin{itemize}
    \item When $\Phi_T$ is ergodic, $\cZ^0(\Phi_T)$ and $\cZ^{\infty}(\Phi_T)$ have the same identity component isomorphic to $\RR\times \SO(n-1)$;
    \item When $\Phi_T$ is non-ergodic, $\cZ^{\infty}(\Phi)\subseteq \cZ^{\infty}(\Phi_T)$ contains a continuously immersed copy of $(C^{\infty}_c((0,1),\RR),+)$.
\end{itemize}
\end{theorem}

We remark here that in the proof of the non-ergodic case of Theorem \ref{T:frame-flow}, we also establish the existence of a nonconstant $\Phi$-invariant smooth function on $FN$. This is based on a structural result of frame flows (see \cite[Corollary 3.11]{Lef23}).

In general, the same equivalence between ergodicity and symmetry also holds for the $k$-frame flow $(F^kN,\mu_k,\Phi^{(k)}_T)$ where $1\leq k\leq n$. To be precise, for $1\le k\leq n$, consider the right action of $L_k:=\mathrm{SO}(n-k)\leq K$ fixing the first $k$ vectors in each orthonormal frame in $FN$, and define
$$
F^kN:=FN/L_k,\qquad\Phi_t^{(k)}(uL_k):=\Phi_t(u)L_k,\qquad \mu_k:=(\pi_k)_*\mu,
$$
where $\pi_k:FN\to F^{(k)}N$ is the natural projection. There are two additional ingredients in the general proof: the mixing criterion for isometric extensions and the smoothness of invariant eigen-sections (see \cite[Theorem 1(b), Proposition 2.7]{Lef23}).

\subsubsection{Automorphisms of Kummer-type complex surfaces} Let $T=\mathbb C^2/\Gamma$ be a complex torus, and let $G$ be a nontrivial finite subgroup of $\operatorname{Aut}(\CC^2,\Gamma)$; for example, fix an elliptic curve $E=\CC/\Lambda$, let $T=E\times E$,  and take $G$ to be generated by the involution $\sigma(z_1,z_2) = (-z_1, -z_2)$. We always assume that $G$ preserves a nonzero holomorphic $2$-form $\Omega$ on $T$. The quotient $Y:=T/G$ is a complex surface with finitely many isolated quotient singularities, and its minimal resolution $X$ is a K3 surface. Denote the quotient map and the resolution map by
$$
q\colon T\to Y,\quad\hbox{ and }\quad \pi\colon X\to Y,
$$
respectively. Let $\mu_X$ denote the normalized smooth volume measure induced by the extended holomorphic $2$-form on $X$, and let $J_0$ denote the standard complex structure on $X$.

Let $f\in\Aut(\CC^2,\Gamma)$ normalize $G$. It induces a holomorphic automorphism $F$ of $X$. The holomorphic automorphism group of a K3 surface is discrete, and
hence so is the holomorphic centralizer of any of its automorphisms \cite[Chapter~15]{Huybrechts}.
For the hyperbolic Kummer examples, a stronger conclusion holds: the cyclic group generated by the automorphism has finite index in its holomorphic centralizer \cite[Proposition~8.1]{Cantat2011}.
The following theorem concerns the full smooth centralizer; as in Theorem \ref{T:main}, we build the bridge between ergodicity and symmetry in this context.

\begin{theorem}\label{T:c-surface}
    For the above system $(X,\mu_X,F)$, the followings are equivalent: \begin{enumerate}
        \item $f$ is hyperbolic;
        \item $F$ is ergodic with respect to $\mu_X$;
        \item $\cZ^{\infty}(F)$ is a finite-dimensional Lie group;
        \item the only essentially bounded $F$-invariant almost complex structures on $X$ are $\pm J_0$.
    \end{enumerate}
Moreover, \begin{itemize}
    \item When $F$ is ergodic, it is mixing, and $\cZ^{\Lip}(F)=\cZ^{\infty}(F)\subseteq \Aut^{\pm}(X)$ is a discrete group;
    \item When $F$ is non-ergodic, $\cZ^{\infty}(F)$ contains a continuously immersed copy of $(C^{\infty}_c((0,1),\RR),+)$.
\end{itemize}
\end{theorem}

This result suggests studying smooth centralizers of automorphisms of other smooth projective varieties, including non-Kummer K3 surfaces; see \cite{CantatSurvey} for background on surface dynamics in particular.
Which dynamical hypotheses ensure that the smooth centralizer is a Lie group, and when must elements of the smooth centralizer be holomorphic?
The connection with invariant complex structures also appears in Zhang's work \cite{Zhang26}, which develops a common cocycle approach to affine centralizer rigidity and rigidity of invariant complex structures.


\subsection{Organization of paper and strategy of proof}

Throughout the paper, we always assume that the subgroup $B$  contains no nontrivial normal subgroups of $G$. This doesn't affect the conclusion since we may change the presentation $(G,B)$ of the homogeneous space $X$ by quotienting out the largest normal subgroup of $G$ contained in $B$.


In Section \ref{S:exist}, we introduce a modification of Dani's suspension of affine maps on compact homogeneous spaces and then give a complete proof of Theorem \ref{T:aux-main}. 
This technique also enables us to prove, as a by-product, that the ergodicity of an affine map is equivalent to the ergodicity of its induced maps on the maximal semisimple and solvable quotient.

In Section \ref{S:proofA}, we present two different approaches to constructing an infinite-dimensional centralizer for a non-ergodic affine map $f=L_a\circ \Phi$ on $X=G/B$.

In Section \ref{SS:key-prop}, we establish a general proposition (Proposition \ref{P:function-vector-field}) for such constructions on connected compact smooth manifolds, whose input consists of a nonconstant invariant smooth function and an invariant smooth vector field with a suitable common nonvanishing set.

In Section \ref{SS:Lie-approach}, we provide a proof with a more Lie-theoretic flavor, which is based on the observation that at least one of the following two possibilities occurs: \begin{enumerate}
    \item[(a)] $\Ad_a\circ D\Phi$ has eigenvalue $1$;
    \item[(b)] $G$ is solvable, and for the nilradical $N$ of $G$, one has
\[\text{either }G=NB,\text{ or the projection of }f \text{ onto }G/NB \text{ is non-ergodic}.
\]
\end{enumerate} 

In Case (a), we start with a $1$-eigenvector of $\Ad_a\circ D\Phi$ to construct an $f$-invariant smooth vector field on $X$, and then produce an abundance of $f$-invariant smooth vector fields by multiplying it by any nonconstant $f$-invariant smooth function. 
Each such $f$-invariant smooth vector field generates a flow commuting with $f$.

In Case (b), we take advantage of the special structure of the solvmanifold $X$, and produce a wild centralizer by using the tower collapse and the two key propositions in \cite{DWWX}. This depends on the classical result that any non-ergodic affine map on a torus admits a periodic quotient torus.

In Section \ref{SS:dyn-approach}, we provide another proof with a more dynamical flavor, which relies on a delicate argument that any nonconstant $f$-invariant smooth function on $X$ yields an $f$-invariant smooth vector field on $X$. This argument makes heavy use of the special structure of affine maps on homogeneous spaces.

In Section \ref{S:bi-homo}, we first prove the non-ergodic case of Theorem \ref{T:bi-homo} by following the arguments in the proof of Theorem \ref{T:main}. The ergodic case of Theorem \ref{T:bi-homo} is obtained by generalizing the arguments in \cite[Proofs of Proposition 6 and Theorem B]{DWWX}, where a new general proposition is established for the Lipschitz centralizer of an affine-induced diffeomorphism to be smooth and finite-dimensional. 

In Section \ref{S:fram-flows}, we first prove the non-ergodic case of Theorem \ref{T:frame-flow} based on a structural result of frame flows in \cite{Lef23}, and then prove the ergodic case of Theorem \ref{T:frame-flow} using the expansiveness of the geodesic flow on $SN$.

In Section \ref{S:Kummer}, we first prove the non-ergodic case of Theorem \ref{T:c-surface} by applying Proposition \ref{P:function-vector-field} again, and then prove the ergodic case of Theorem \ref{T:c-surface} based on some regularity arguments in \cite[Section 3]{Zhang26}.

\subsection{Acknowledgement}

We are thankful to Simion Filip for helpful discussions around the proof of Lemma \ref{L:inv-vec-field}. DD was supported by Swedish Research Council grant VR 2023-03596.
AW was supported by NSF DMS-2452616. DX was supported by 2024YFA 1015100, NSFC 12090010 and 12090015.

\subsection{Declaration of AI use}

This paper was written by humans, for humans.  After more than a year of human effort devoted to understanding separate cases, the main result in this paper (Theorem~\ref{T:main}) was finally proved in early September, 2026, without AI assistance.  In the process of checking our original proof, an AI assistant found a significant way to shorten it (this shortened proof is presented in Section~\ref{SS:Lie-approach}). The same AI assistant found a second proof in mid-October (the second proof is presented in Section~\ref{SS:dyn-approach}). AI subsequently helped develop the non-affine examples presented in Sections~\ref{S:fram-flows} and \ref{S:Kummer} (although our proofs deviate significantly from the AI suggestions). 

The authors proudly take responsibility for the math in this paper.

\section{Existence of a nonconstant invariant smooth function}\label{S:exist}

\subsection{Modifying Dani's suspension}\label{SS:Danisus}
In this subsection, we introduce a modification of Dani's construction of suspensions of affine diffeomorphisms on homogeneous spaces. This useful tool allows us to transfer certain properties of homogeneous translations and flows to the setting of affine diffeomorphisms. The original source of these ideas can be found in \cite[Section 7]{Dani77}, although the construction there requires correction.

Let $G$ be a connected, simply connected Lie group, and let $B$ a closed subgroup of $G$ such that $X=G/B$ admits a finite $G$-invariant measure. Let $f=L_a\circ \Phi$ be an affine diffeomorphism of $X$, where $a\in G$, $\Phi\in\Aut(G)$ and $\Phi(B)=B$. The normalized $G$-invariant measure on $X$, denoted by $m_X$, is preserved by $f$.

Since $G$ is simply connected, $\Aut(G)\cong\Aut(\fg)$ is an algebraic group. Write $c_a$ for the conjugation map by $a$. Let $\cA$ denote the algebraic hull of $c_a\circ \Phi$ in $\Aut(G)$. Then $\cA$ is abelian and has finitely many connected components. Choose a positive integer $m$ such that $(c_a\circ \Phi)^m\in\cA^0$, and choose a one-parameter subgroup $\{\Phi_t:t\in\RR\}$ of $\cA^0$ such that $(c_a\circ \Phi)^m=\Phi_1$.

We consider the semi-direct product $\widetilde{G}:=G\rtimes\RR$, where the multiplication law is given by $$
    (g,t)\cdot(h,s):=(g\cdot\Phi_t(h),t+s),
    \quad \forall g,h\in G,\;t,s\in\RR.
$$
Write $a_m=a\Phi(a)\cdots \Phi^{m-1}(a)$ and $\delta:=(a_m^{-1},1)\in \widetilde{G}$. Define $\widetilde{B}$ to be the subgroup generated by $B\times\{0\}$ and $\delta$. It is easy to check that $\delta(B\times \{0\})\delta^{-1}=B\times\{0\}$, which implies that $\widetilde{B}$ is a closed subgroup of $\widetilde{G}$. Write $\widetilde{X}=\widetilde{G}/\widetilde{B}$. Note that there is an almost everywhere identification 
$$
    X\times [0,1)\xrightarrow{\cong} \widetilde{X},\quad (gB,t)\mapsto (\Phi_t(g),t)\widetilde{B};
    $$ 
The Haar measure on $\widetilde{X}$ equals the push forward of the product measure on $X\times [0,1)$.

The group homomorphism $G\to\widetilde{G}$, $g\mapsto(g,0)$, induces a smooth embedding $$
\iota\colon X\to\widetilde{X},\quad gB\mapsto(g,0)\widetilde{B}.
$$
Moreover, for $\widetilde{g}:=(1_G,1)\in\widetilde{G}$, we see that $L_{\widetilde{g}}\circ \iota=\iota\circ f^m$ on $X$. This completes our modified construction of suspension of affine maps. 

Finally, we briefly explain two reasons for introducing this modification of Dani's construction. First, Dani's original construction realizes $f$ as the restriction of a homogeneous translation on a larger space, but it does not in general provide a one-parameter subgroup containing the translation. In our construction, after passing to a positive iterate, $f^m$ is realized as the time-one map of a homogeneous flow on $\widetilde{X}$, restricted to the embedded copy $\iota(X)$. This allows us to use results for homogeneous flows in the proof of Theorem \ref{T:aux-main}. Second, the proof of \cite[Theorem 7.1]{Dani77} invokes \cite[Lemma 1.10]{Dani77}, which is actually false as stated. 
Our proof of the dichotomy for the ergodicity of affine maps avoids invoking this lemma.

\subsection{Proof of Theorem \ref{T:aux-main}}
In this subsection, we give a complete proof of Theorem \ref{T:aux-main}, which concerns the the density of invariant smooth functions in the space of all invariant $L^2$ functions for affine maps on homogeneous spaces. 

We start with a general lemma about the density of invariant smooth functions in the space of all invariant $L^2$ functions for homogeneous flows.

\begin{lemma}\label{L:smoothdense}
Let $S$ be a connected Lie group, and $D$ a closed subgroup of $S$ such that $Z=S/D$ admits a finite $S$-invariant measure. Let $F=\{a_t:t\in\RR\}$ be a one-parameter subgroup of $S$. Then the subspace $C^{\infty}(Z)^F\cap L^2(Z)$ of all $F$-invariant smooth $L^2$ functions on $Z$ is dense in the space $L^2(Z)^F$ of all $F$-invariant $L^2$ functions on $Z$.
\end{lemma}
\begin{proof}
By lifting $S$ to its universal cover and replacing $D$ by its full preimage, we may assume that $S$ is simply connected. Let $H$ denote the Moore normal subgroup associated with $F=\{a_t:t\in\RR\}$, that is, the smallest closed connected normal subgroup of $S$ such that
\begin{equation}\label{E:relcpt}
\{\Ad(a_t):t\in\RR\}\text{ is relatively compact in }
\Aut(\Lie(S/H)).
\end{equation}
It follows from Moore-Mautner's theorem (see \cite[Theorem 1.1]{Moore80}) that
$$
    L^2(Z)^F\subseteq L^2(Z)^H\cong L^2(S/\overline{HD}).
$$

By \eqref{E:relcpt}, we may choose an inner product on $\Lie(S/H)$ that is $\Ad(a_t)$-invariant for every $t\in\RR$. This inner product uniquely determines a right-invariant Riemannian metric on the Lie group $S/H$, and hence induces a Riemannian metric on the homogeneous space
$$
    \overline{Z}:=S/\overline{HD}\cong\frac{S/H}{\overline{HD}/H}.
$$
Under this metric, the left translation by $a_t\,(t\in\RR)$ is an isometry on $\overline{Z}$. Let $K$ denote the closure of $\{L_{a_t}:t\in\RR\}$ in $\Isom(\overline{Z})$. It is easy to see that $K$ is compact.

Let $\Psi\in L^2(Z)^F$ and write its projection onto $\overline{Z}$ as $\overline{\Psi}\in L^2(\overline{Z})$. We choose a sequence of smooth compactly supported functions $\{\psi_n:n\in\NN\}$ on $\overline{Z}$ such that $\psi_n\to\overline{\Psi}$ in $L^2(\overline{Z})$, and write
$$
    \phi_n(z):=\int_K\psi_n(kz)\mathrm{d}k,
    \quad z\in \overline{Z},
$$
where $\mathrm{d}k$ is the normalized Haar measure on $K$. Then each $\phi_n$ is an $F$-invariant smooth function on $\overline{Z}$, and $\phi_n\to \overline{\Psi}$ in $L^2(\overline{Z})$. Pulling back to $Z$ gives the desired $L^2$ density.
\end{proof}

\begin{remark}
    In fact, the subspace $C^{\omega}(Z)^F\cap L^2(Z)$ of all $F$-invariant analytic $L^2$ functions on $Z$ is also dense in $L^2(Z)^F$. This can be done by choosing a sequence of analytic $L^2$ functions $\{\psi_n:n\in\NN\}$ on $\overline{Z}$ such that $\psi_n\to\overline{\Psi}$ in $L^2(\overline{Z})$ (see \cite[Lemma 42.7]{KM97}). 
\end{remark}


\begin{proposition}\label{P:smoothdense-affine-finite-volume}
Let $G$ be a connected Lie group, and $B$ a closed subgroup of $G$ such that $X=G/B$ admits a finite $G$-invariant measure. Let $f=L_a\circ\Phi$ be an affine diffeomorphism of $X$, where $a\in G$, $\Phi\in\Aut(G)$ and $\Phi(B)=B$. Then the subspace $C^{\infty}(X)^f\cap L^2(X)$ is dense in $L^2(X)^f$.

In particular, $f$ is non-ergodic if and only if there exists a nonconstant $f$-invariant smooth $L^2$ function on $X$.
\end{proposition}
\begin{proof} By lifting $G$ to its universal cover, replacing $B$ by its full preimage, and lifting $a$ and $\Phi$, we may assume that $G$ is simply connected. We keep the notations for modified Dani's suspension of the affine map $f$ on $X$ as in Section \ref{SS:Danisus}. 

Let $\Psi\in L^2(X)^f$. We define $$
    \widetilde{\Psi}\colon \widetilde{X}\to \CC,\quad (g,t)\widetilde{B}\mapsto \Psi(\Phi_{-t}(g)B).
    $$
    It follows from our construction that $\widetilde{\Psi}$ is a well-defined $L^2$ function on $\widetilde{X}$, and it is invariant under the one-parameter subgroup $F=\{\widetilde{g}_s=(1_G,s):s\in\RR\}$. Then Lemma \ref{L:smoothdense} gives a sequence $\{\psi_n:n\in\NN\}$ of $F$-invariant smooth $L^2$ functions on $\widetilde{X}$ such that $\psi_n\to \widetilde{\Psi}$ in $L^2(\widetilde{X})$. By the $F$-invariances of $\psi_n$ and $\widetilde{\Psi}$, restricting to the zero section $X$ gives $\psi_n\circ \iota\to \Psi$ in $L^2(X)$. In particular, $$
    \frac{1}{m}\sum_{i=0}^{m-1}\psi_n\circ \iota\circ f^i\to \Psi \text{ in }L^2(X).
    $$
    Note that the left-hand side is a sequence of $f$-invariant smooth $L^2$ functions on $X$. 
    
    Finally, if $f$ is non-ergodic, we may choose $\Psi$ to be nonconstant. Then in the above approximation, when $n$ is large enough, the function $\frac{1}{m}\sum_{i=0}^{m-1}\psi_n\circ \iota\circ f^i$ is also nonconstant. This completes the proof.
\end{proof}

For later use, we also prove here the density of smooth $L^2$ eigenfunctions in the space of all $L^2$ eigenfunctions for affine maps on homogeneous spaces.

\begin{corollary}\label{C:dense-eigen-smooth}
Let $(G,B,X,f)$ be as in Proposition \ref{P:smoothdense-affine-finite-volume}, and let $U_f$ denote the Koopman operator associated with $f$ on $L^2(X)$. Then for any $\zeta\in\CC$ with $|\zeta|=1$, the subspace of all smooth $L^2$ $\zeta$-eigenfunctions for $U_f$ is dense in the space of all $L^2$ $\zeta$-eigenfunctions for $U_f$.
\end{corollary}
\begin{proof}
Let $\Psi\in L^2(X)$ satisfy $\Psi\circ f=\zeta\cdot \Psi$. Consider the homogeneous space $$
X\times S^1=(G\times S^1)/(B\times \{1\})
$$
and the affine diffeomorphism $$
\widehat{f}\colon X\times S^1\to X\times S^1,\qquad (x,z)\mapsto (f(x),\zeta^{-1}\cdot z).
$$
It follows that the $L^2$ function $$
\widehat{\Psi}\colon X\times S^1\to \CC,\qquad (x,z)\mapsto\Psi(x)\cdot z
$$
is $\widehat{f}$-invariant. Then Proposition \ref{P:smoothdense-affine-finite-volume} gives a sequence $\{\widehat{\psi}_n:n\in\NN\}$ in $C^{\infty}(X\times S^1)^{\widehat{f}}\cap L^2(X\times S^1)$ with $\|\widehat{\psi}_n-\widehat{\Psi}\|_{L^2(X\times S^1)}\to 0$ as $n\to +\infty$. For $n\in\NN$, write $$
\psi_n\colon X\to \CC,\qquad x\mapsto \int_{S^1}\widehat{\psi}_n(x,z)\cdot\overline{z}\mathrm{d}z.
$$
It follows that each $\psi_n$ is smooth and $\psi_n\circ f=\zeta\cdot\psi_n$. Moreover, $$
\|\psi_n-\Psi\|_{L^2(X)}\leq \|\widehat{\psi}_n-\widehat{\Psi}\|_{L^2(X\times S^1)}\to 0\quad \text{ as }n\to +\infty.$$
This completes the proof.  
\end{proof}

\subsection{Dichotomy for ergodicity of affine maps}
As an application of our modified Dani's suspension, we establish the dichotomy for ergodicity on the maximal semisimple and solvable quotients. Since our suspension is a suspension of a power of the affine map, we first prove a total ergodicity statement.

\begin{lemma}\label{L:affine-total-ergodic-finite-volume}
Let $S$ be a connected Lie group, let $D$ be a closed subgroup such that $Y=S/D$ admits a finite $S$-invariant measure, and let $f$ be an ergodic affine diffeomorphism of $Y$. Then $f^q$ is ergodic for all $q\in\NN$.
\end{lemma}
\begin{proof}
Suppose that $f^q$ is not ergodic for some $q\in\NN$. This means that $L^2(Y)^{f^q}\neq \CC$. Write $\cH:=L^2(Y)^{f^q}$. It is clear that the restriction of the Koopman operator $U_f$ on $\cH$ has finite order. Since $f$ is ergodic, we see that $U_f\big|_{\cH}\neq \id_{\cH}$. It follows that there exists a nonzero $\zeta$-eigenfunction for $U_f$, where $\zeta^q=1$ and $\zeta\ne1$. By Corollary \ref{C:dense-eigen-smooth}, we may take this eigenfunction $u$ to be smooth. Then $u^q$ is $f$-invariant, and hence is constant almost everywhere by ergodicity. By continuity, we see that $u$ is a nonzero constant, which contradicts $u\circ f=\zeta u$ and $\zeta\ne1$.
\end{proof}

Let $G$ be a connected Lie group, let $B$ be a closed subgroup such that $X=G/B$ admits a finite $G$-invariant measure, and let $R$ be the radical of $G$. Put $G^{(0)}=G$, $G^{(j+1)}=[G^{(j)},G^{(j)}]$, and let $G^{(\infty)}$ be the stable term of this derived series. We denote the maximal semisimple and solvable quotients of $X$ by
$$
X_{ss}:=G/\overline{RB},\qquad
X_{sol}:=G/\overline{G^{(\infty)}B}.
$$
Both subgroups $R$ and $G^{(\infty)}$ are characteristic, so every affine diffeomorphism of $X$ induces affine diffeomorphisms on $X_{ss}$ and $X_{sol}$.

\begin{proposition}\label{P:affine-ergodicity-dichotomy-A}
Let $f=L_a\circ\Phi$ be an affine diffeomorphism of $X$ as above. Then $f$ is ergodic if and only if its induced maps on both $X_{ss}$ and $X_{sol}$ are ergodic.
\end{proposition}
\begin{proof}
The ``only if'' part is clear. For the ``if'' part, by passing to the universal cover of $G$ and replacing $B$ by its full inverse image, we may assume that $G$ is simply connected. Let $\widetilde{G}$, $\widetilde{B}$, $\widetilde{X}$ and $F=\{(1_G,s):s\in\RR\}$ be the suspension of $f^m$ from Section \ref{SS:Danisus}. It carries a finite $\widetilde{G}$-invariant measure, and its flow is ergodic if and only if $f^m$ is ergodic on $X$.

For $N=R$ or $G^{(\infty)}$, since $N$ is characteristic in $G$, it is normal in $\widetilde{G}$. Moreover, we have $\overline{N\widetilde{B}}\cap G=\overline{NB}.$ Put
$$
\widetilde{X}_N:=\widetilde{G}/\overline{N\widetilde{B}}.
$$
Then $\widetilde{X}_N$ is the suspension of the map induced by $f^m$ on $G/\overline{NB}$. In particular, its $F$-flow is ergodic if and only if that induced map is ergodic. By assumption, the maps induced by $f$ on $X_{ss}$ and $X_{sol}$ are ergodic. By Lemma \ref{L:affine-total-ergodic-finite-volume}, their $m$-th powers are ergodic, and hence that the flow $F$ is ergodic on both $\widetilde{X}_R$ and $\widetilde{X}_{G^{(\infty)}}$.

On one hand, let $\widetilde{R}$ be the radical of $\widetilde{G}$. Since $R\subseteq\widetilde{R}$, the maximal semisimple quotient $\widetilde{G}/\overline{\widetilde{R}\widetilde{B}}$ is a factor of $\widetilde{X}_R$, so the flow $F$ on it is ergodic.

On the other hand, we first show that $G^{(\infty)}=\widetilde{G}^{(\infty)}$. Indeed, it is clear that $G^{(\infty)}\subseteq \widetilde{G}^{(\infty)}$. Conversely, since $G^{(\infty)}$ is normal in $\widetilde{G}$, and $\widetilde{G}/G^{(\infty)}$ is an abelian extension of the solvable group $G/G^{(\infty)}$, we see that $\widetilde{G}/G^{(\infty)}$ is solvable, so $\widetilde{G}^{(\infty)}\subseteq G^{(\infty)}$. It follows that the maximal solvable quotient of $\widetilde{X}$ is precisely $\widetilde{X}_{G^{(\infty)}}$, on which the flow is also ergodic.

The ergodicity criterion for homogeneous flows on finite-volume spaces (see \cite[Corollary 2]{Starkov87}) now implies that $\widetilde{F}$ is ergodic on $\widetilde{X}$. Therefore, $f^m$ is ergodic on $X$, and by Lemma \ref{L:affine-total-ergodic-finite-volume}, so is $f$. The proof is complete.
\end{proof}

\section{Two proofs of Theorem \ref{T:main}}\label{S:proofA}

\subsection{A key proposition}\label{SS:key-prop}
In this subsection, we establish a general proposition for constructing infinite-dimensional centralizers of non-ergodic diffeomorphisms. 
The proof of this key proposition depends on the newly proven existence of a nonconstant invariant smooth function (see Theorem \ref{T:aux-main}).


\begin{proposition}\label{P:function-vector-field}
    Let $M$ be a connected compact smooth manifold, and $f$ a diffeomorphism on $M$. Suppose that there exists a nonconstant $f$-invariant smooth function $\Psi$ on $M$, an $f$-invariant smooth vector field $\fX$ on $M$, and a small open neighborhood $\cU\subseteq M$ in which every point $p$ satisfies $$
\fX_p\neq 0\quad\text{and}\quad (\mathrm{d}\Psi)_p\neq 0.
$$ Then there is a continuous injective group homomorphism $(C_c^{\infty}((0,1),\RR),+)\hookrightarrow \cZ^{\infty}(f).$
    In particular, $\cZ^{\infty}(f)$ isn't a finite-dimensional Lie group.
\end{proposition}
\begin{proof}
We divide the proof into several steps.
\medskip

\textsl{Step 1: More $f$-invariant vector fields with disjoint supports.}
By shrinking $\cU$, we may choose a nonempty open interval $J\subseteq \RR$ such that
$$\overline J\subset\Psi(\cU)\cap\mathrm{int}(\Psi(M)).$$
We choose pairwise disjoint nonempty open subintervals $J_n\subset J\,(n\geq 1)$, and choose nonzero real functions $\chi_n\in C_c^\infty(J_n,\RR)$ extended by zero to $\mathbb R$. Define
\[
 \mathfrak X^{(n)}:=\chi_n(\Psi)\cdot \mathfrak X.
\]
We see that each $\mathfrak X^{(n)}$ is nonzero. Indeed, choose $s\in J_n$ with $\chi_n(s)\ne0$. Since $s\in\Psi(\cU)$, there is an $x\in \cU$ with $\Psi(x)=s$, and then $\mathfrak X^{(n)}_x=\chi_n(s)\cdot \fX_x\neq 0$.

For each $n\geq 1$, we see that for all $p\in M$, $$
Df_p(\fX^{(n)}_p)=\chi_n(\Psi(p))\cdot Df_p(\fX_p)=\chi_n(\Psi(f(p)))\cdot \fX_{f(p)}=\fX^{(n)}_{f(p)},
$$
so $\fX^{(n)}$ is $f$-invariant. Since $\supp(\fX^{(n)})\subseteq \Psi^{-1}(\supp(\chi_n))$, the supports of the vector fields $\fX^{(n)}\,(n\geq 1)$ are pairwise disjoint. It follows that
\begin{equation*}\label{eq:disjoint-brackets}
 [\fX^{(n)},\fX^{(m)}]=0
 \qquad(n,m\ge1).
\end{equation*}
So these vector fields $\fX^{(n)}\,(n\geq 1)$ generate nontrivial commuting flows $\{\phi^{(n)}_t:t\in\RR\}$ on $X$, and all these flows commute with $f$.\vspace{1em}

\textsl{Step 2: Exclusion of nonzero global periods.}
We claim that for each $n\geq 1$, \begin{equation}
    \phi^{(n)}_t=\id_M\Longrightarrow t=0.
\end{equation}
In fact, set $\cO:=M\setminus \Psi^{-1}(\overline{J})$. This is a nonempty open subset. Let $n\geq 1$. Since $\fX^{(n)}$ vanishes on $\cO$, we see that for any $t\in\RR$, the map $\phi^{(n)}_t$ fixes all points in $\cO$. Suppose that $\phi^{(n)}_T=\id_M$ for some $T\neq 0$. Without loss of generality, we may assume that $T>0$. So this flow defines a smooth action $$
(\RR/T\ZZ)\times M\to M,\quad (t,p)\mapsto \phi^{(n)}_t(p).
$$
We choose any smooth Riemannian metric $g_0$ on $M$, and define its average as $$
\overline{g_0}:=\frac{1}{T}\int_0^T(\phi^{(n)}_t)^*g_0\mathrm{d}t.
$$
This is again a smooth Riemannian metric on $M$, under which each $\phi^{(n)}_t\,(t\in\RR)$ is an isometry. However, an isometry fixing a nonempty open set pointwise must be the identity on the connected space. This forces that the flow $\{\phi^{(n)}_t:t\in\RR\}$ is trivial, which is a contradiction.\vspace{1em}

\textsl{Step 3: Construction of the additive group homomorphism.}
Fix a Riemannian metric and a connection on $TM$, and use the increasing norms of vector fields on $M$:
\[
 \|\fX\|_{C^m}:=\max_{0\le j\le m}\sup_{p\in M}|\nabla^j\fX_p|, \qquad m\ge0.
\]
For $n\geq 1$, choose
\[
 \varepsilon_n=\frac{2^{-n}}{1+\|\mathfrak X^{(n)}\|_{C^n}}.
\]
Let $(t_n)_{n\ge1}$ be a dense sequence in $(0,1)$. For $\rho\in C_c^\infty((0,1),\mathbb R)$, define
\begin{equation}\label{eq:series-field}
 \mathfrak Y^{\rho}:=\sum_{n=1}^{\infty}\varepsilon_n \rho(t_n)\cdot \mathfrak X^{(n)}.
\end{equation}
For every fixed $m\ge0$ and $n\ge\max\{1,m\}$, one has
\[
 \|\varepsilon_n \rho(t_n)\cdot \mathfrak{X}^{(n)}\|_{C^m}
 \le\|\rho\|_{\infty}2^{-n}.
\]
Thus the series \eqref{eq:series-field} converges in every $C^m$ norm and defines a smooth vector field on all of $M$. Moreover, for each $m\geq 0$ there is a finite constant $\kappa_m$ independent of $\rho$ such that
\begin{equation}\label{eq:continuity-estimate}
 \|\mathfrak Y^{\rho}\|_{C^m}\le \kappa_m\|\rho\|_{\infty}.
\end{equation}
In particular, $\rho\mapsto\mathfrak Y^{\rho}$ is linear and continuous for the usual test-function topology.

The continuity of pushforward by $f$ gives
$f_*\mathfrak Y^{\rho}=\mathfrak Y^{\rho}$. By
\eqref{eq:disjoint-brackets} and continuity of the Lie bracket
under smooth convergence, we obtain that
\begin{equation}\label{eq:commuting-series}
 [\mathfrak Y^{\rho_1},\mathfrak Y^{\rho_2}]=0,
 \quad\text{for all }\rho_1,\rho_2\in C_c^\infty((0,1),\mathbb R).
\end{equation}
All these vector fields are complete because $M$ is compact. Write $\{\varphi^{\rho}_t:t\in\RR\}$ for the flow generated by $\fY^{\rho}$, and define
\[
 \iota(\rho):=\varphi^{\rho}_1.
\]
Then each $\iota(\rho)$ commutes with $f$, and $\iota(\rho_1)\circ \iota(\rho_2)=\iota(\rho_1+\rho_2)$ for all $\rho_1,\rho_2$ in $C_c^\infty((0,1),\mathbb R)$. It follows that $\iota\colon C_c^\infty((0,1),\mathbb R)\to \cZ^{\infty}(f)$ is a group homomorphism. Estimate \eqref{eq:continuity-estimate}, together with continuous dependence of flows on smooth vector fields on a compact manifold, proves that $\iota$ is continuous in the stated topologies.\vspace{1em}

\textsl{Step 4: Injectivity of $\iota$.}
Suppose that $\iota(\rho)=\id_M$ for some $\rho\in C_c^\infty((0,1),\mathbb R)$. Let $n\geq 1$. Then the disjointness of the supports gives
\[
 \mathfrak Y^{\rho}=\varepsilon_n \rho(t_n)\cdot \mathfrak X^{(n)}\quad \text{on}\quad \supp(\fX^{(n)}).
\]
Since the flow $\{\phi^{(n)}_t:t\in\RR\}$ preserves $\supp(\fX^{(n)})$, the uniqueness of integral curves gives
\[
 \varphi^{\rho}_1(p)=\phi^{(n)}_{\varepsilon_n\rho(t_n)}(p)
 \qquad\text{for all }p\in \supp(\fX^{(n)}).
\]
Since $\varphi^{\rho}_1=\id_M$, the right-hand map is the identity on $\supp(\fX^{(n)})$. It is also the identity outside that support,
and so
\[
 \phi^{(n)}_{\varepsilon_n\rho(t_n)}=\id_M.
\]
Since $\varepsilon_n>0$, Step 2 then implies that $\rho(t_n)=0$. This holds for every $n$. Then the density of $(t_n)_{n\in\NN}$ in $(0,1)$ and continuity of $\rho$ imply that $\rho\equiv0$. Hence $\iota$ is injective.
\end{proof}


\subsection{A Lie-theoretical proof of Theorem \ref{T:main}}\label{SS:Lie-approach}
In this section, we present our first approach to constructing infinite-dimensional centralizers for non-ergodic affine maps on compact homogeneous spaces. This approach has a more Lie-theoretic flavor and consists of the following two parts.

Let $f=L_a\circ \Phi$ be a non-ergodic affine diffeomorphism on $X=G/B$. Then at least one of the following two possibilities holds (see Lemma \ref{L:dichotomy}):
\begin{enumerate}
    \item[(a)] $\Ad_a\circ D\Phi$ has eigenvalue $1$;
    \item[(b)] $G$ is solvable, and for the nilradical $N$ of $G$, one has
$$\text{either }G=NB,\text{ or the projection of }f \text{ onto }G/NB \text{ is non-ergodic}.
$$
\end{enumerate} 
We prove Theorem \ref{T:main} by using different methods in these two possibilities.

\subsubsection{An elementary Lie-theoretical lemma} We prove the above alternatives as:

\begin{lemma}\label{L:dichotomy}
    Let $G$ be a connected Lie group, and $B$ a closed cocompact subgroup of $G$ such that $X=G/B$ admits a finite $G$-invariant measure. Suppose that the subgroup $B$ contains no nontrivial connected normal subgroup of $G$.

    Let $f$ be a non-ergodic affine diffeomorphism of $X$. Then one of the two possibilities (a) and (b) holds.
\end{lemma}

To prove Lemma \ref{L:dichotomy}, we first prove an elementary result about automorphisms of solvable groups.

\begin{lemma}\label{L:torbaselinear}
Let $R$ be a connected solvable Lie group, $N$ the nilradical of $R$, and $\Phi\in\Aut(R)$. Then the induced automorphism
$\overline{\Phi}\in\Aut(R/N)$ has finite order.
\end{lemma}
\begin{proof}
Write $\fr=\Lie(R)$ and $\fn=\Lie(N)$, and complexify these Lie algebras as $\fr_{\CC}=\fr\otimes_{\RR}\CC$ and $\fn_{\CC}=\fn\otimes_{\RR}\CC$. By Lie's theorem, there is a basis of $\fr_{\CC}$ in which all matrices $\ad(X)\,(X\in\fr_{\CC})$, are upper triangular. Write $m=\dim R$ and let $\alpha_1,\ldots,\alpha_m$ denote their diagonal entries, regarded as complex-linear functionals on $\fr_{\CC}$.

We note that
$$\fn_{\CC}=\bigcap_{i=1}^m\ker(\alpha_i).$$
Consequently, the induced functionals $\overline{\alpha_i}$ on $\fr_{\CC}/\fn_{\CC}$ satisfy
\begin{equation}\label{E:weightlinearspan}
\Span_{\CC}\{\overline{\alpha_i}:1\leq i\leq m\}=(\fr_{\CC}/\fn_{\CC})^*.
\end{equation}
On the other hand, for every $X\in\fr_{\CC}$, one has
$$\ad(D\Phi_{\CC}(X))=D\Phi_{\CC}\circ\ad(X)\circ D\Phi_{\CC}^{-1}.$$
So the two adjoint operators have the same characteristic polynomial, that is,
$$\prod_{i=1}^m\bigl(t-\alpha_i(D\Phi_{\CC}(X))\bigr)=\prod_{i=1}^m\bigl(t-\alpha_i(X)\bigr).
$$
This is an identity of polynomials in $t$ and the coordinates of $X$. Unique factorization shows that $D\Phi_{\CC}^{\top}$ permutes the multiset $\{\alpha_1,\ldots,\alpha_m\}$. Thus, there exists $r\geq1$ such that
$$\alpha_i\circ(D\Phi_{\CC})^r=\alpha_i,
\qquad 1\leq i\leq m.
$$
By \eqref{E:weightlinearspan}, the automorphism induced by $(D\Phi_{\CC})^r$ on $\fr_{\CC}/\fn_{\CC}$ is identity. It follows that $\overline{D\Phi}^r=\id$ on $\fr/\fn$.

Since $R/N$ is connected, an automorphism of $R/N$ is determined by its differential. Therefore, we conclude that $\overline{\Phi}^{\,r}=\id$.
\end{proof}

\begin{corollary}\label{C:trivial-linear}
Let $R$ be a connected solvable Lie group, $N$ its nilradical, and $B\leq R$ a closed cocompact subgroup such that $X=R/B$ admits a finite $R$-invariant measure. Suppose that the subgroup $B$ doesn't contain any nontrivial connected normal subgroup of $R$.

Let $f=L_a\circ \Phi$ be an affine diffeomorphism of $X$ that projects to an ergodic map on the nontrivial space $R/NB$. Then the induced automorphism of $\Ad_a\circ D\Phi$ on $\Lie(R/N)$ is identity. 
\end{corollary}
\begin{proof}
It is clear that the induced automorphism of $\Ad_a$ on $\Lie(R/N)$ is identity. Moreover, Lemma \ref{L:torbaselinear} says that the induced automorphism of $D\Phi$ on $\Lie(R/N)$ has finite order. We need to show that it is identity.

In fact, since $R/NB$ is a nontrivial connected compact abelian group, we may identify it with $\TT^k$ for some positive integer $k$, and write the induced map of $f$ as $$
\overline{f}\colon \TT^k\to\TT^k,\quad v\mapsto\alpha+Uv,
$$
where $\alpha\in\RR^k$ and $U\in\GL(k,\ZZ)$. The preceding paragraph gives $U^r=\id$ for some $r\geq 1$. It follows that $\overline f^{\,r}$ is the translation by
$$
(\id+U+\cdots+U^{r-1})\alpha\in\ker(U-\id).
$$
Define a subtorus
\[
\TT_0:=\bigl(\ker(U-\id)+\ZZ^k\bigr)/\ZZ^k.
\]
Since $U(\TT_0)=\TT_0$, the map $\overline f$ induces an
affine map on $\TT^k/\TT_0$, and the $r$-th power of this
map is identity. If $U\neq \id$, then $\TT^k/\TT_0$ is a positive-dimensional torus, so the induced periodic map $\overline{f}$ cannot be ergodic. 
Hence $U=\id$.
\end{proof}

Now we proceed to the proof of Lemma \ref{L:dichotomy}:

\begin{proof}[Proof of Lemma \ref{L:dichotomy}] Let $R$ denote the solvable radical of $G$. If $G$ is not solvable, then $G/R$ is a nontrivial connected semisimple Lie group. It follows from Jacobson's fixed-point theorem \cite[Theorem~9]{Jacobson} that the induced automorphism of $\Ad_a\circ D\Phi$ on $\Lie(G/R)$ has eigenvalue $1$. So the automorphism $\Ad_a\circ D\Phi$ on $\Lie(G)$ also has eigenvalue $1$, that is, the possibility (a) holds.

If $G=R$ is solvable and $\Ad_a\circ D\Phi$ doesn't have eigenvalue $1$, then by Corollary \ref{C:trivial-linear}, we see that either $R=NB$, or the projection of $f$ onto $R/NB$  is non-ergodic. This means that the possibility (b) holds.
\end{proof}

\subsubsection{Proof in the possibility (a)} In this paragraph, we prove Theorem \ref{T:main} in the possibility (a), where we assume that the Lie algebra automorphism $\Ad_a\circ D\Phi$ has eigenvalue $1$. The following proof is a direct application of Proposition \ref{P:function-vector-field}.



\begin{proof}[Proof of Theorem \ref{T:main} in the possibility (a)] Since $f$ is non-ergodic, Theorem \ref{T:aux-main} gives a nonconstant $f$-invariant smooth function $\Psi$ on $X$. By Proposition \ref{P:function-vector-field}, it suffices to construct an $f$-invariant smooth vector field $\fX$ on $X$ whose vanishing set has no interior.

Let $v\in \Lie(G)$ be a nonzero vector such that $(\Ad_a\circ D\Phi)v=v$. We define a vector field on $X$ by $$
\fX_{gB}:=\left.\frac{\mathrm{d}}{\mathrm{d}t}\right|_{t=0}\exp(tv)gB,\quad \forall g\in G.
$$
We first note that $$
f(\exp(tv)gB)=\exp(t(\Ad_a\circ D\Phi)v)\cdot f(gB)=\exp(tv)\cdot f(gB),\quad\forall g\in G,\;t\in\RR.
$$
Taking the derivative at $t=0$ gives $Df(\fX_{gB})=\fX_{f(gB)}$, that is, $\fX$ is $f$-invariant. We also observe that $\fX$ is a nonzero vector field. Otherwise, one would have $$
\exp(\RR v)\subseteq\bigcap_{g\in G}gBg^{-1},
$$
whose right-hand side is a discrete group by our assumption. This is a contradiction. So the vector field $\fX$ generates a nontrivial flow on $X$ commuting with $f$.

Moreover, we claim that the set of nonvanishing points of $\fX$ is open dense in $X$. In fact, for $g\in G$, one has $$
\fX_{gB}=0\Longleftrightarrow \Ad_{g^{-1}}(v)\in \Lie(B).
$$
If the vector field $\fX$ vanishes on a nonempty open subset of $X$, then the map $$
G\to \Lie(G)/\Lie(B),\quad g\mapsto \Ad_{g^{-1}}(v)+\Lie(B)
$$
vanishes on a nonempty open subset of $G$. Since $G$ is connected and the above map is real analytic, this forces the above map to vanish on all of $G$. This implies that $\fX\equiv 0$, which is a contradiction. This completes the proof.
\end{proof}

\subsubsection{Proof in the possibility (b)}
 In this paragraph, we prove Theorem \ref{T:main} in the possibility (b), where we assume that $G=R$ is solvable, and for the nilradical $N$ of $R$, one has
$$\text{either }R=NB,\text{ or the projection of }f \text{ onto }R/NB \text{ is non-ergodic}.
$$
 We always assume that the subgroup $B$ doesn't contain any nontrivial connected normal subgroup of $R$. Here we note that $B^0\subseteq N$ and that $NB$ is closed in $R$ (see \cite[Section 5]{Mos54}). 


The proof here is based on the two key propositions in \cite{DWWX}. We start with a standard result for a non-ergodic affine map on a torus:

\begin{lemma}\label{L:torus}
    Let $f\colon x\mapsto Ax+b$ be a non-ergodic affine map on $\TT^k$, where $A\in\GL(k,\ZZ)$ and $b\in\RR^k$. Then there exists a proper $A$-invariant subtorus $H^0\subseteq \TT^k$ such that the projection of $f$ onto $\TT^k/H^0$ is periodic.
\end{lemma}
\begin{proof}
Since $f$ is non-ergodic, there is a nonconstant function $u\in L^2(\TT^k)$ such that $u\circ f=u$. Consider the Fourier expansion of $u$:
\[
u=\sum_{m\in\mathbb Z^k}c_m\chi_m,\qquad\text{where }
\chi_m(x)=e^{2\pi i\langle m,x\rangle}.
\]
Since $\chi_m\circ f=\chi_m(b)\cdot \chi_{A^{\mathsf T}m}$, comparing Fourier coefficients gives
\begin{equation}\label{E:Fourier}
    c_{A^{\mathsf T}m}=\chi_m(b)\cdot c_m.
\end{equation}
Choose $m\in\ZZ^k\setminus\{0\}$ such that $c_m\neq0$. All coefficients along its $A^{\mathsf T}$-orbit have the same nonzero absolute value. Since $\sum_{m\in\mathbb Z^k}|c_m|^2<\infty$,
this orbit is finite. Thus there exists $r\geq1$ such that
$(A^{\mathsf T})^r m=m$.

Put
\[
m_j=(A^{\mathsf T})^jm\text{ for }j=0,1,\cdots,r-1,\quad\text{and}\quad
b_r=b+Ab+\cdots+A^{r-1}b.
\]
Multiplying \eqref{E:Fourier} over $m=m_0,m_1,\cdots,m_{r-1}$ gives
\[
c_m
=\left(\prod_{j=0}^{r-1}\chi_{m_j}(b)\right)c_m
=\chi_m(b_r)\cdot c_m.
\]
Hence $\chi_m(b_r)=1$. Reading the indices modulo $r$, we also have
\[
\chi_{m_j}(b_r)
=\prod_{i=0}^{r-1}\chi_{m_{j+i}}(b)
=\prod_{i=0}^{r-1}\chi_{m_i}(b)
=1,
\qquad j=0,\ldots,r-1.
\]
Since $f^r(x)=A^rx+b_r$ and $(A^{\mathsf T})^rm_j=m_j$, it follows that
\begin{equation}\label{E:character}
    \chi_{m_j}\circ f^r=\chi_{m_j},\qquad j=0,\ldots,r-1.
\end{equation}

Define
\[H:=\bigcap_{j=0}^{r-1}\ker\chi_{m_j}.\]
Then $H^0$ is a proper $A$-invariant subtorus, and hence $f$ induces an affine map $\overline f$ on the positive-dimensional $\TT^k/H^0$. 
We shall show that it is periodic. In fact, \eqref{E:character} gives
\[f^r(x)-x\in H,\quad\forall  x\in \TT^k.\]
Consequently, the following continuous map taking discrete values
\[\TT^k\to H/H^0,\qquad x\mapsto f^r(x)-x+H^0\]
must be constant.
Denote its value by $\eta\in H/H^0$. Then it follows that
\[
\overline f^{\,r}(y)=y+\eta,\quad\forall y\in \TT^k/H^0.
\]
For $k=[H:H^0]$, we have $k\eta=0$, and hence $\overline f^{kr}=\operatorname{id}_{\TT^k/H^0}.$  This is as desired.
\end{proof}

Now we proceed to the proof of Theorem \ref{T:main} in the possibility (b).

\begin{proof}[Proof of Theorem \ref{T:main} in the possibility (b)]\,

\textsl{Case 1: $R=NB$.}
 In this case, since $B^0\subseteq N$, we see that $R=N$ is nilpotent. Now $B$ is a closed cocompact subgroup of $N$. In view of \cite[Corollary 4.3]{Witte87}, we know that $B^0$ is normal in $N$, and hence is trivial by our assumption. So the subgroup $B$ is a cocompact lattice in $N$. 

Since the affine map $f$ on $X=N/B$ is non-ergodic, it follows from \cite[Theorem 1]{Pa69} that its projection on $\overline{X}:=N/[N,N]B$ is also non-ergodic. We identify the positive-dimensional torus $\overline{X}$ with $\TT^k$ for some positive integer $k$.

By Lemma \ref{L:torus}, there is a tower of $f$-invariant fiber bundles 
\begin{equation}
\begin{tikzcd}
  {[N,N]B/B} \arrow[r] & X=N/B \arrow[d] \\
  H^0\arrow[r] & \overline{X}=N/[N,N]B\cong \TT^k \arrow[d] \\
  & \TT^k/H^0
\end{tikzcd}
\end{equation}
and $f$ projects to a periodic map on the base. In view of \cite[Corollary 29]{DWWX}, $\cZ^{\infty}(f)$ contains a subgroup isomorphic to either $C_c^{\infty}((0,1),\RR)$ or $\Diff_c^{\infty}((0,1))$. In particular, $\cZ^{\infty}(f)$ isn't a finite-dimensional Lie group.
\medskip

\textsl{Case 2: the projection of $f$ onto $R/NB$ is non-ergodic.} In this case, we have the $f$-invariant Mostow's sequence (see \cite[Section 5]{Mos54}) $$
NB/B\to R/B\to R/NB,
$$
and the projection of $f$ onto the base torus $R/NB$ is non-ergodic. The remaining arguments using Lemma \ref{L:torus} and \cite[Corollary 29]{DWWX} are the same as in Case 1.
\end{proof}

\subsection{A dynamical proof of Theorem \ref{T:main}}\label{SS:dyn-approach}
In this section, we provide another approach to proving  Theorem \ref{T:main}, which has its own dynamical features compared to our first approach. 
In this proof, we can even obtain the following stronger result, which is the ``only if'' part of Theorem \ref{T:main-prime}. (Note that the ``if'' part of Theorem \ref{T:main-prime} is proved in \cite[Theorem B]{DWWX}.)

\begin{proposition}\label{P:key-2}
    Let $G$ be a connected Lie group, and $B$ a closed cocompact subgroup of $G$ such that $X=G/B$ admits a finite $G$-invariant measure. Suppose that the subgroup $B$ contains no nontrivial connected normal subgroup of $G$.

    Let $f=L_a\circ \Phi$ be a non-ergodic affine diffeomorphism of $X$. Then there is a continuous injective group homomorphism $(C_c^{\infty}((0,1),\RR),+)\hookrightarrow \cZ^{\infty}(f).$
    In particular, $\cZ^{\infty}(f)$ isn't a finite-dimensional Lie group.
\end{proposition}

\subsubsection{Existence of a nonzero invariant smooth vector field}

To prove Proposition \ref{P:key-2}, in view of Theorem \ref{T:aux-main} and Proposition \ref{P:function-vector-field}, it suffices to prove that: a nonconstant $f$-invariant smooth function $\Psi$ on $X$ yields an $f$-invariant smooth vector field $\fX$ on $X$ with a suitable common nonvanishing set. We isolate this result as follows:

\begin{lemma}\label{L:inv-vec-field}
    Let $G,B,X,f$ be as in Proposition \ref{P:key-2}. Suppose that there exists a nonconstant $f$-invariant smooth function $\Psi$ on $X$. Then there exists an $f$-invariant smooth vector field $\fX$ on $X$ and a small open neighborhood $\cU\subseteq X$ in which every point $x$ satisfies $$
\fX_x\neq 0\quad\text{and}\quad (\mathrm{d}\Psi)_x\neq 0.
$$
\end{lemma}
\begin{proof}
Write $\fg=\Lie(G)$ and $A=\Ad_a\circ D\Phi\in \Aut(\fg)$. For any $x\in X$, define $$
j_x\colon \fg\to T_xX,\qquad v\mapsto \left.\frac{\mathrm{d}}{\mathrm{d}t}\right|_{t=0}\exp(tv)x.
$$
Each $j_x$ is a surjective linear map. It is straightforward to see that $j_{f(x)}\circ A=Df_x\circ j_x$. Moreover, differentiating $\Psi\circ f=\Psi$ gives $\mathrm{d}\Psi_{f(x)}\circ Df_x=\mathrm{d}\Psi_x$. Define a smooth map $$
b\colon X\to \fg^*,\qquad x\mapsto \mathrm{d}\Psi_x\circ j_x.
$$
It follows that \begin{equation}\label{E:b-f-A}
    b(f(x))\circ A=\mathrm{d}\Psi_{f(x)}\circ j_{f(x)}\circ A=\mathrm{d}\Psi_{f(x)}\circ Df_x\circ j_x=\mathrm{d}\Psi_x\circ j_x=b(x).
\end{equation}

Now we fix a basis $\{e_1,\cdots,e_n\}$ of $\fg$, and write $b_j(x)=b(x)e_j$ for $j=1,\cdots,n$. We claim that: the finite-dimensional subspace of $C^{\infty}(X,\CC)$ $$
E:=\Span_{\CC}\{b_1(\cdot),\cdots,b_n(\cdot)\}
$$
is invariant under the Koopman operator associated with $f$. In fact, write $A^{-1}e_j=\sum_{i=1}^nA^{ij}e_i$. Then for $j=1,\cdots,n$, one has $$
b_j(f(x))=b(f(x))e_j\overset{\eqref{E:b-f-A}}{=}(b(x)\circ A^{-1})e_j=\sum_{i=1}^nA^{ij}\cdot b(x)e_i=\sum_{i=1}^nA^{ij}\cdot b_i(x).
$$
This verifies the claim.

It follows that the invariant subspace $E$ admits a basis consisting of eigenfunctions of the unitary Koopman operator associated with $f$, say $\{u_1(\cdot),\cdots,u_k(\cdot)\}$, where \begin{equation}\label{E:u-eigen}
u_i(\cdot)\in C^{\infty}(X,\CC),\qquad u_i\circ f=\zeta_i u_i,\quad\zeta_i\in\CC \text{ with } |\zeta_i|=1.
\end{equation}
Write $b_j(\cdot)=\sum_{i=1}^kc_{ij}u_i(\cdot)$ for $c_{ij}\in\CC$, and $\ell_i:=(c_{i1},\cdots,c_{in})\in \CC^{1\times n}\cong \fg_{\CC}^*$. It is straightforward to see that \begin{equation}\label{E:b-u-relation}
b(x)=\sum_{i=1}^ku_i(x)\ell_i.
\end{equation}
Then we have \begin{equation*}
    \begin{aligned}
        0\overset{\eqref{E:b-f-A}}{=}b(f(x))\circ A-b(x)&\overset{\eqref{E:b-u-relation}}{=}\sum_{i=1}^k(u_i(f(x))\cdot \ell_i\circ A-u_i(x)\ell_i)\\
        &\overset{\eqref{E:u-eigen}}{=}\sum_{i=1}^ku_i(x)\cdot (\zeta_i\cdot\ell_i\circ A-\ell_i).
    \end{aligned}
\end{equation*}
By the linear independence of $\{u_1(\cdot),\cdots,u_k(\cdot)\}$, one has $$\ell_i\circ A=\zeta_i^{-1}\ell_i \quad \text{for}\quad i=1,\cdots,k.$$

We choose and fix $x_0\in X$ with $(\mathrm{d}\Psi)_{x_0}\neq 0$. Then we see that $b(x_0)\neq 0$, and hence that there exists some $1\leq i\leq k$ such that $u_i(x_0)\ell_i\neq 0$ in \eqref{E:b-u-relation}. In particular, $u_i(x_0)\neq 0$ and $\ell_i\neq 0$. This implies that $\zeta_i^{-1}$ is an eigenvalue of $A$. Let $v_0\in\fg_{\CC}$ be a nonzero vector with $Av_0=\zeta_i^{-1}v_0$. Define $$
\fX_x=\overline{u_i(x)}\cdot j_x(v_0),\qquad x\in X.
$$
Then $\fX$ is a smooth complex vector field on $X$, and it is $f$-invariant: \begin{equation*}
    \begin{aligned}
\fX_{f(x)}=\overline{u_i(f(x))}\cdot j_{f(x)}(v_0)&=\overline{\zeta_iu_i(x)}\cdot Df_x\circ j_x\circ A^{-1}(v_0)\\
&=\overline{\zeta_i}\cdot Df_x(\overline{u_i(x)}\cdot j_x(\zeta_iv_0))\\
&=Df_x(\fX_x).
\end{aligned}
\end{equation*}

Now consider the nonempty open set $$
\cU_0:=\left\{x\in X: (\mathrm{d}\Psi)_{x}\neq 0\text{ and }u_i(x)\neq 0\right\}.
$$
Since the set $\{x\in X:j_x(v_0)\neq 0\}$ is open and dense in $X$ (see Proof of Theorem \ref{T:main} in the possibility (a)), it contains some point $x_0'\in \cU_0$. It follows that $\fX_{x_0'}\neq 0$. Taking the real or imaginary part of $\fX$ nonvanishing at $x_0'$ gives the desired one.
\end{proof}

\begin{remark}
    It is clear that Lemma \ref{L:inv-vec-field} holds for isometries of closed Riemannian manifolds and symplectomorphisms of closed symplectic manifolds, by taking, respectively, the gradient and the Hamiltonian vector field of a nonconstant smooth invariant function.

    However, it doesn't hold for a general smooth diffeomorphism, even if it preserves a smooth volume. In the next paragraph, we will present a simple counterexample of a partially hyperbolic, volume-preserving diffeomorphism on $\TT^3$ that has a nonconstant invariant smooth function but doesn't admit a nonzero invariant smooth vector field.
\end{remark}

\subsubsection{Two examples where Lemma \ref{L:inv-vec-field} holds or fails}

The first example corresponds to \cite[Example 22]{DWWX} where a non-ergodic affine map on $\TT^3$ has an Anosov fiber map. We explicitly construct its invariant smooth functions and invariant smooth vector fields.

\begin{example}
Let $X=\TT^3$, and write $$
f\colon \TT^3\to\TT^3,\qquad (t,x,y)\mapsto (t+1/2, t+2x+y,x+y).
$$
It is clear that the function $$
\Psi\colon \TT^3\to \RR,\qquad (t,x,y)\mapsto \cos(4\pi t).
$$
is a nonconstant, $f$-invariant, smooth function on $\TT^3$. Moreover, the vector field $$
\fX=\frac{\partial}{\partial t}-\frac{\partial}{\partial y}
$$
is a nonzero, $f$-invariant, smooth vector field on $\TT^3$.
\end{example}

The second example is a special non-affine perturbation of a non-ergodic affine map on $\TT^3$. We show that it admits a nonzero invariant smooth function, but doesn't admit any nonzero invariant vector field.

\begin{example}
Let $X=\TT^3$, and write $$
f_0\colon \TT^3\to\TT^3,\qquad (t,x,y)\mapsto (t,2x+y,x+y).
$$
For $\epsilon>0$, we consider the following small perturbation: $$
L_{\epsilon}\colon \TT^3\to\TT^3,\qquad (t,x,y)\mapsto \left(t,\frac{\epsilon}{2\pi}\cos(2\pi t)\sin(2\pi y)+x,y\right).
$$
When $\epsilon>0$ is sufficiently small, the composition $f_{\epsilon}:=L_{\epsilon}\circ f_0$ is a partially hyperbolic, volume-preserving diffeomorphism on $\TT^3$. It is clear that the function
$$
\Psi\colon \TT^3\to\RR,\qquad (t,x,y)\mapsto \cos(2\pi t).
$$
is a nonconstant, $f_{\epsilon}$-invariant, smooth function on $\TT^3$.

Now we claim that: the map $f_{\epsilon}$ doesn't admit any nonzero invariant smooth vector field on $\TT^3$. Suppose that $\fX$ is an $f_{\epsilon}$-invariant smooth vector field on $\TT^3$. In the coordinate $(t,z)=(t,(x,y))\in\TT^3$, write \begin{equation*}
        \fX_{(t,z)}=\left(a(t,z),\fY_{(t,z)}\right),\qquad \text{where }\fY\text{ is a smooth vector field on }\TT^2.
    \end{equation*}
We also write $f_{\epsilon}(t,z)=(t,g_{\epsilon,t}(z))$, where $$
        g_{\epsilon,t}(x,y):=\left(
        \frac{\epsilon}{2\pi}\cos(2\pi t)\sin(2\pi(x+y))+(2x+y),\;x+y
        \right).
 $$
Using $Df_{(t,z)}\fX_{(t,z)}=\fX_{f(t,z)}$, we see that \begin{subequations}
    \begin{align}
        a(t,z)&=a(t,g_{\epsilon,t}(z));\label{E:a-inv}\\
        a(t,z)\cdot\frac{\partial g_{\epsilon,t}(z)}{\partial t}+D_zg_{\epsilon,t}(z)\cdot \fY_{(t,z)}&=\fY_{(t,\,g_{\epsilon,t}(z))}.\label{E:Y-inv}
    \end{align}
\end{subequations}

On one hand, when $\epsilon>0$ is sufficiently small, we know that for all $t\in\TT$, $g_{\epsilon,t}(\cdot)$ is a transitive Anosov diffeomorphism on $\TT^2$. Then it follows from \eqref{E:a-inv} that the function $a(t,z)$ is independent of $z$. In particular, one has $a(t,z)\equiv a(t,0)$ for all $(t,z)\in\TT^3$.

On the other hand, taking $z=0$ in \eqref{E:Y-inv} gives $$
D_zg_{\epsilon,t}(z)\big|_{z=0}\cdot \fY_{(t,0)}=\fY_{(t,0)}.
$$
Note that $$B_t:=D_zg_{\epsilon,t}(z)\big|_{z=0}=\begin{pmatrix}
    \epsilon\cos(2\pi t)+2 & \epsilon\cos(2\pi t)+1\\ 1 &1
\end{pmatrix}.$$
When $0<\epsilon<1$, we know that for all $t\in\TT$ this matrix $B_t$ has no eigenvalue $1$, which forces that $\fY_{(t,0)}\equiv 0$. 

Moreover, differentiating at $z=0$ in \eqref{E:Y-inv} gives $$
a(t,0)\cdot\frac{\partial B_t}{\partial t}+B_t\cdot (D_z\fY)(t,0)=(D_z\fY)(t,0)\cdot B_t.
$$
Then taking trace gives $$
0=a(t,0)\cdot\frac{\partial \tr(B_t)}{\partial t}=a(t,0)\cdot (-2\pi\epsilon\sin(2\pi t)).
$$
Since $-2\pi\epsilon\sin(2\pi t)\neq 0$ at all but finitely many points of $t\in\TT$, by continuity, we see that $a(t,0)\equiv 0$. Then it follows from \eqref{E:Y-inv} that $$
D_zg_{\epsilon,t}(z)\cdot \fY_{(t,z)}=\fY_{(t,\,g_{\epsilon,t}(z))}.
$$
Since an Anosov diffeomorphism admits no nonzero continuous invariant vector field, we conclude that $\fY_{(t,z)}\equiv 0$. This proves that $\fX_{(t,z)}\equiv 0$.
\end{example}

\section{Application (I): Affine-induced maps on bi-homogeneous spaces}\label{S:bi-homo}

In this section, we prove Theorem \ref{T:bi-homo} as an analogue of Theorem \ref{T:main} in the setting of affine-induced diffeomorphisms on bi-homogeneous spaces.

We first fix some notations. Let $X=G/B$ be a compact homogeneous space with a $G$-invariant probability measure, and let $K\leq G$ be a compact subgroup acting freely by left translations on $X$. Write $\fg=\Lie(G)$, $\fk=\Lie(K)$. Let $f=L_a\circ \Phi$ be an affine diffeomorphism on $X$ with $(c_a\circ \Phi)(K)=K$. It induces a diffeomorphism $F$ on $Y=K\backslash G/B$. Write $\pi\colon X\to Y$ to be the natural projection map. Then we have $F\circ \pi=\pi\circ f$. 

\subsection{The non-ergodic case} To prove the non-ergodic case of Theorem \ref{T:bi-homo}, in view of Proposition \ref{P:function-vector-field}, it suffices to show the existence of the followings: \begin{enumerate}
    \item[(i)] a nonconstant $F$-invariant smooth function $\Psi$ on $Y$;
    \item[(ii)] an $F$-invariant invariant smooth vector field $\fX$ on $Y$, and a small open neighborhood $\cU\subseteq Y$ in which every point $y$ satisfies $$
\fX_y\neq 0\quad\text{and}\quad (\mathrm{d}\Psi)_y\neq 0.
$$
\end{enumerate}
We complete these two tasks in two lemmas below.

\begin{lemma}\label{L:bi-homo-inv}
Let $F$ be an affine-induced diffeomorphism on $Y=K\backslash G/B$. Then the subspace $C^{\infty}(Y)^F$ of all $F$-invariant smooth functions on $Y$ is dense in the space $L^2(Y)^F$ of all $F$-invariant $L^2$ functions on $Y$.

    In particular, $F$ is non-ergodic if and only if there exists a nonconstant $F$-invariant smooth function on $Y$.
\end{lemma}
\begin{proof}
Let $\Psi\in L^2(Y)^F$. Then $\Psi\circ \pi\in L^2(X)^f$. By Theorem \ref{T:aux-main}, we choose a sequence of $f$-invariant smooth functions $\{\psi_n:n\in\NN\}$ on $X$ such that $\psi_n\to\Psi\circ \pi$ in $L^2(X)$, and write
$$
    \overline{\psi}_n(x):=\int_K\psi_n(kx)\mathrm{d}k,
    \quad x\in X,
$$
where $\mathrm{d}k$ is the normalized Haar measure on $K$. We claim that: each $\overline{\psi}_n$ is also an $f$-invariant smooth function on $X$. In fact, for $f=L_a\circ \Phi\in \Aff(G,B)$ with $(c_a\circ \Phi)(K)=K$, write $\theta:=(c_a\circ \Phi)\big|_K\in\Aut(K)$. It follows that for any $x\in X$, $$
\overline{\psi}_n(f(x))=\int_K\psi_n(kf(x))\mathrm{d}k=\int_K(\psi_n\circ f)(\theta^{-1}(k)x)\mathrm{d}k=\int_K\psi_n(kx)\mathrm{d}k=\overline{\psi}_n(x).
$$

By construction, each $\overline{\psi}_n$ descends to an $F$-invariant smooth function on $Y$ as $\phi_n$. Then we see that $\|\phi_n-\Psi\|_{L^2(Y)}\leq \|\psi_n-\Psi\circ \pi\|_{L^2(X)}\to 0$ as $n\to +\infty$. This completes the proof.
\end{proof}

\begin{lemma}
    Let $G,B,K,Y,F$ be as in Lemma \ref{L:bi-homo-inv}. Suppose that there exists a nonconstant $F$-invariant smooth function $\Psi$ on $Y$. Then there exists an $F$-invariant smooth vector field $\fY$ on $Y$ and a small open neighborhood $\cU\subseteq Y$ in which every point $y$ satisfies $$
\fY_y\neq 0\quad\text{and}\quad (\mathrm{d}\Psi)_y\neq 0.
$$
\end{lemma}
\begin{proof} We divide the proof into several steps.
\medskip

\textsl{Step 1: the smooth map induced by the differential.}

Write $A=(\Ad_a\circ D\Phi)\big|_{\fg/\fk}\in \Aut(\fg/\fk)$. For any $x\in X$, define $$
j_x\colon \fg/\fk\to T_{\pi(x)}Y,\qquad v+\fk\mapsto \left.\frac{\mathrm{d}}{\mathrm{d}t}\right|_{t=0}\pi(\exp(tv)x).
$$
Each $j_x$ is a surjective linear map. It is straightforward to see that $$
j_{kx}\circ \Ad_k\big|_{\fg/\fk}=j_x\;(k\in K);
\qquad
j_{f(x)}\circ A=DF_{\pi(x)}\circ j_x.
$$ Moreover, differentiating $\Psi\circ F=\Psi$ gives $\mathrm{d}\Psi_{F(y)}\circ DF_y=\mathrm{d}\Psi_y$ for any $y\in Y$. Define a smooth map $$
b\colon X\to (\fg/\fk)^*,\qquad x\mapsto \mathrm{d}\Psi_{\pi(x)}\circ j_x.
$$
It follows that \begin{equation}\label{E:bx-prop}
    b(kx)\circ \Ad_k\big|_{\fg/\fk}=b(x)\;(k\in K);\qquad
    b(f(x))\circ A=b(x).
\end{equation}
Moreover, we have $A\circ \Ad_k\circ A^{-1}=\Ad_{(c_a\circ \Phi)(k)}$ for all $k\in K$.\medskip 

\textsl{Step 2: the real multiplicative Jordan decomposition.}

Consider the real multiplicative Jordan decomposition of $A\in \GL(\fg/\fk)$ (see \cite[Section 4.3]{Morris05} for details): $$
A=A_e\circ A_h\circ A_u,
$$
where $A_e,A_h,A_u\in \GL(\fg/\fk)$ are the elliptic, hyperbolic, unipotent part of $A$, and they commute with each other. Note that $A$ normalizes the compact subgroup 
$$
C:=\left\{\Ad_k\big|_{\fg/\fk}:k\in K\right\}\leq \GL(\fg/\fk).
$$
We first claim that: $A_e$ normalizes $C$, and $A_h,A_u$ centralizes $C$. In fact, consider the adjoint representation of $\GL(\fg/\fk)$ on the vector space $\End(\fg/\fk)$. Then $$
\Ad_A=\Ad_{A_e}\circ \Ad_{A_h}\circ \Ad_{A_u}\text{ preserves the subspace }S:=\Span(C)\subseteq \End(\fg/\fk).
$$
Moreover, since $\Ad_A^n(C)\subseteq C$ for any $n\in\ZZ$, we see that $\Ad_{A_h}\big|_S=\Ad_{A_u}\big|_S=\id_S$. This verifies the claim.

We write $$
V_0:=\ker(A_h-\id),\qquad N:=\log\left(A_u\big|_{V_0}\right).
$$
It follows from the above claim that the subspace $V_0$ is invariant under $A$ and $C$, and that the nilpotent operator $N$ commutes with $A_e\big|_{V_0}$ and $C\big|_{V_0}$. We next claim that: for any $x\in X$, \begin{enumerate}
    \item[(b1)] the linear functional $b(x)$ vanishes on all eigenspaces of $A_h$ except $V_0$;
    \item[(b2)] $b(x)\big|_{V_0}\circ N=0$.
\end{enumerate}

To verify (b1), we consider the eigenspace decomposition for the action $A_h$ on $\fg/\fk$. Let $v\in \fg/\fk$ be a $\lambda$-eigenvector of $A_h$ where $\lambda>0$ and $\lambda\neq 1$. Without loss of generality, we assume that $0<\lambda<1$. Then we have $A^n(v)\to 0$ as $n\to +\infty$. It follows from the second identity in \eqref{E:bx-prop} that $$
|b(x)(v)|=|b(f^n(x))\circ A^n(v)|\leq 
\sup_{x'\in X}\|b(x')\|\cdot \|A^n(v)\|\to 0\text{ as }n\to +\infty.
$$
This gives $b(x)v=0$. 

To verify (b2), we note that $A\big|_{V_0}=A_e\big|_{V_0}\circ A_h\big|_{V_0}\circ A_u\big|_{V_0}=A_e\big|_{V_0}\circ\exp(N)$. It follows from the second identity in \eqref{E:bx-prop} that for any $n\in\ZZ$, $$
b(f^n(x))\circ A^n_e\big|_{V_0}=b(f^n(x))\circ A^n\big|_{V_0}\circ \exp(-nN)=b(x)\big|_{V_0}\circ \exp(-nN).
$$
The left-hand side has bounded operator norms in $n$, while the right-hand side has polynomial operator norms in $n$. So we conclude that $b(x)\big|_{V_0}\circ N=0$. 

Therefore, the second identity in \eqref{E:bx-prop} on $V_0$ becomes \begin{equation}\label{E:bx-prop-2}
\begin{aligned}
b(f(x))\circ A_e\big|_{V_0}&=b(f(x))\circ A\big|_{V_0}\circ A_h^{-1}\big|_{V_0}\circ\exp(-N)\\
&=b(x)\big|_{V_0}\circ \exp(-N)=b(x)\big|_{V_0}.
\end{aligned}
\end{equation}
\medskip

\textsl{Step 3: the construction of invariant smooth vector fields.} 

Since both $C$ and $A_e$ commutes with $A_h$, the compact subgroup $$
D:=\overline{\langle C,A_e\rangle}\leq \GL(\fg/\fk)
$$
also preserves the eigenspace $V_0=\ker(A_h-\id)$. We choose and fix any $D$-invariant inner product $\langle\cdot,\cdot\rangle$ on $V_0$. For any $x\in X$, let $w(x)\in V_0$ be the unique vector with $$
\langle w(x),v\rangle=b(x)v,\qquad \forall v\in V_0.
$$  
It is straightforward to check that $$
w(kx)=\Ad_k(w(x))\;(k\in K);\qquad w(f(x))=A_e(w(x)).
$$
For $r\in\NN_0$ and $x\in X$, define $$\fX^{(r)}_x:=\left(j_x\big|_{V_0}\circ N^r\right)(w(x))\in T_{\pi(x)}Y.$$ 
Note that for any $k\in K$, $\fX^{(r)}_{kx}=\fX^{(r)}_{x}$. It follows that the smooth correspondence $x\mapsto \fX^{(r)}_x$ descends to a smooth vector field on $Y=\pi(X)$, denoted by $\fY^{(r)}$. 

We first observe that $\fY^{(0)}\not\equiv 0$. In fact, for $x\in X$, $$
\mathrm{d}\Psi_{\pi(x)}(\fY^{(0)}_{\pi(x)})=\mathrm{d}\Psi_{\pi(x)}(\fX^{(0)}_x)=\mathrm{d}\Psi_{\pi(x)}\circ j_x(w(x))=b(x)(w(x))=\|w(x)\|^2.
$$
Then we see that $\mathrm{d}\Psi_{\pi(x)}=0\Longleftrightarrow b(x)=0\Longleftrightarrow w(x)=0\Longleftrightarrow \mathrm{d}\Psi_{\pi(x)}(\fY^{(0)}_{\pi(x)})=0$. In particular, whenever $\mathrm{d}\Psi_y\neq 0$, one has $\fY^{(0)}_y\neq 0$.

Write $r_0=\max\{r\geq 0:\fY^{(r)}\not\equiv 0\}$. We claim that: $\fY^{(r_0)}$ is a nonzero $F$-invariant smooth vector field on $Y$. In fact, for any $r\in\NN_0$ and $x\in X$, \begin{equation*}
    \begin{aligned}
        DF_{\pi(x)}\circ \fY^{(r)}_{\pi(x)}&=
        \left(DF_{\pi(x)}\circ j_x\big|_{V_0}\circ N^r\right)(w(x))\\
        &=\left(j_{f(x)}\circ A\big|_{V_0}\circ N^r\right)(w(x))\\
        &=\left(j_{f(x)}\circ A_u\big|_{V_0}\circ N^r\right)(A_e(w(x)))\\
        &=\left(j_{f(x)}\circ \exp(N)\circ N^r\right)(w(f(x)))\\
        &=\sum_{\ell\geq 0}\frac{\fX^{(r+\ell)}_{f(x)}}{\ell!}
        =\sum_{\ell\geq 0}\frac{\fY^{(r+\ell)}_{F(\pi(x))}}{\ell!}.
    \end{aligned}
\end{equation*}
By the maximality of $r_0$, we have $DF_{\pi(x)}\circ \fY^{(r_0)}_{\pi(x)}=\fY^{(r_0)}_{F(\pi(x))}$. This verifies the claim. Moreover, whenever $\fY^{(r_0)}_{\pi(x)}\neq 0$, one has $w(x)\neq 0$, that is, $\mathrm{d}\Psi_{\pi(x)}\neq 0$.
\end{proof}

\subsection{The ergodic case} 

We first provide a sufficient condition for the Lipschitz centralizer of the ergodic affine-induced diffeomorphism to be smooth and finite-dimensional. To prove this, we start with a general regularity result.

\begin{lemma}\label{L:regularity}
    Let $F$ be an ergodic diffeomorphism on $Y=K\backslash G/B$ induced by an affine map $f=L_a\circ \Phi$ on $X=G/B$, and write $\theta:=(c_a\circ \Phi)\big|_K\in\Aut(K)$. Let $T$ be an invertible linear map on a finite-dimensional vector space $W$. Suppose that there exists a smooth representation $\rho\colon K\to\GL(W)$ with \begin{equation}\label{E:k-w}
    T\circ \rho(k)=\rho(\theta(k))\circ T,\qquad \text{for all }k\in K.
    \end{equation}
    Then the space of all bounded measurable maps $\cA\colon X\to W$ with \begin{subequations}
    \begin{align}
        \cA_{f(x)}&=T(\cA_x),\qquad \text{for }\mu_{\Haar}\text{-a.e. }x\in X;\label{E:a-t-1}\\
        \cA_{kx}&=\rho(k)(\cA_x),\qquad \text{for all }k\in K\text{ and }\mu_{\Haar}\text{-a.e. }x\in X\label{E:a-t-2}
    \end{align}
    \end{subequations}
    is finite-dimensional, of dimension at most $\dim{W}$, and every member has a smooth representative.
\end{lemma}
\begin{proof} We divide the proof into several steps.\medskip

\textsl{Step 1: An $L^2$ norm on the space of bounded measurable maps.}

We choose and fix any norm $\|\cdot\|$ on $W$ induced by an inner product. Let $$
W_b:=\left\{v\in W:\sup_{n\in\ZZ}\|T^nv\|<+\infty\right\}.
$$
Then $W_b$ is a $T$-invariant subspace of $W$. Moreover, we see that the eigenvalues of $T\big|_{W_b}$ are on the unit circle. By \eqref{E:k-w}, the subspace $W_b$ is also invariant under $\rho(K)$, and hence the subgroup \begin{equation*}
    \cG:=\left\{
    \rho(k)\circ T^n\big|_{W_b}:k\in K,\;n\in\ZZ
    \right\}\leq \GL(W_b)
\end{equation*} 
has compact closure. Then we may assume that the inner product $\langle\cdot,\cdot\rangle$ on $W_b$ is $\cG$-invariant by averaging it over the compact closure of $\cG$. We claim that: \begin{equation}\label{E:L-infty-L_2}
\begin{aligned}
    &\left\{
\cA\in L^{\infty}(X,W):\cA \text{ satisfies \eqref{E:a-t-1} and \eqref{E:a-t-2}}
\right\}\\
=&\left\{
\cA\in L^2(X,W_b):\cA \text{ satisfies \eqref{E:a-t-1} and \eqref{E:a-t-2}}
\right\}.
\end{aligned}
\end{equation}

On one hand, for any $\cA\in L^{\infty}(X,W)$ with \eqref{E:a-t-1} and \eqref{E:a-t-2}, we have $\cA_{f^n(x)}=T^n(\cA_x)$ for all $n\in\ZZ$ and $\mu_{\Haar}$-almost all $x\in X$. It follows from the essential boundedness of $\cA$ that there exists a $\mu_{\Haar}$-full subset of $X$ in which each $x$ satisfies $$
\sup_{n\in\ZZ}\|T^n(\cA_x)\|<+\infty,\quad \text{that is},\quad \cA_x\in W_b.
$$
Then by the compactness of $X$, we have $\cA\in L^2(X,W_b)$.

On the other hand, write $\cS$ to be the right-hand side of \eqref{E:L-infty-L_2}. Then for any $\cA,\cA'\in \cS$, their pointwise scalar product \begin{equation}\label{E:scalar-product}
    \Phi\colon X\to \RR,\qquad x\mapsto \langle\cA_x,\cA_x'\rangle
\end{equation}
is invariant under the $K$-action and the affine map $f$ on $X$. It follows that $\Phi$ descends to an $F$-invariant function on $Y$. By the ergodicity of $F$, we see that $\Phi$ is $m$-almost everywhere constant. In particular, for any $\cA\in \cS$, its norm $x\mapsto \|\cA_x\|$ is $m$-almost everywhere constant. This proves that $\cA\in L^{\infty}(X,W)$.
\medskip 

\textsl{Step 2: The smooth approximation in a finite-dimensional space.}

We first verify that: our desired space $\cS$ is finite-dimensional. In fact, let $\cA^{(1)},\cdots$, $\cA^{(N)}$ be $L^2$-orthonormal elements in $\cS$. Then for any $1\leq i<j\leq N$, one has $$
\delta_{ij}=\langle\cA^{(i)},\cA^{(j)}\rangle\equiv \langle\cA^{(i)}_x,\cA^{(j)}_x\rangle\quad \text{for }m\text{-almost all }x\in X.
$$
It follows that for $m$-almost all $x\in X$, $\cA_x^{(1)},\cdots$, $\cA_x^{(N)}$ are orthonormal vectors in $W_b$. In particular, this implies that $N\leq \dim{W_b}$. Hence $\dim{\cS}\leq \dim{W_b}$.

Next, we verify that: the space $\cS$ has a dense subspace $$
\cS':=\left\{\cA\in C^{\infty}(X,W_b):\cA \text{ satisfies \eqref{E:a-t-1} and \eqref{E:a-t-2}}
\right\}.
$$
In fact, let $\cA\in\cS$. Since $T\big|_{W_b}$ is orthogonal, with respect to an orthonormal basis of $(W_b)_{\CC}$, we see from \eqref{E:a-t-1} that any coordinate $a_x$ of $\cA_x$ satisfies $$
a_{f(x)}=\zeta\cdot a_x,\qquad \text{for }\mu_{\Haar}\text{-a.e. }x\in X\text{ and some }\zeta\in\CC\text{ with }|\zeta|=1.
$$
By Corollary \ref{C:dense-eigen-smooth}, we see that each function $x\mapsto a_x$ can be approximated by smooth $\zeta$-eigenfunctions for the Koopman operator associated with $f$ in $L^2(X)$. Therefore, the map $\cA$ can also be approximated by smooth maps satisfying \eqref{E:a-t-1} in $L^2(X,W_b)$. By averaging them over $K$-actions, we may assume that these smooth maps also satisfy \eqref{E:a-t-2}. This is as desired.

Finally, since a dense subspace of a finite-dimensional space must be the whole space, we conclude that $\cS=\cS'$. This completes the proof.
\end{proof}

\begin{proposition}\label{P:Lip-smooth}
    Let $F$ be an ergodic, affine-induced diffeomorphism on $Y=K\backslash G/B$. Suppose that there exists a finite-dimensional real vector space $E_0^*$ in $\Gamma^{\infty}(\Sym^2T^*Y)$ such that $F^*(E_0^*)=E_0^*$ and $E_0^*$ contains a Riemannian metric. Then $\cZ^{\Lip}(F)=\cZ^{\infty}(F)$ is a finite-dimensional Lie group with $\dim\cZ^{\infty}(F)\leq \dim{Y}$.
\end{proposition}
\begin{proof}
We divide the proof into several steps.\medskip

\textsl{Step 1: Construction of the finite-dimensional Hom space.}

Write $A=(\Ad_a\circ D\Phi)\big|_{\fg/\fk}\in \Aut(\fg/\fk)$. Consider the finite-dimensional vector space $$
W:=\Hom(E_0^*,\Sym^2((\fg/\fk)^*)),
$$
and the invertible linear map 
$T\colon W\to W$ defined by \[T(L):= (\Sym^2A^{-1})^*\circ L\circ F^*.\]
There is a smooth representation $$
\rho\colon K\to \GL(W),\qquad k\mapsto (\Sym^2\Ad_k^{-1})^*.
$$
Then it is readily checked that \eqref{E:k-w} holds. 

For any $x\in X$, define $$
j_x\colon \fg/\fk\to T_{\pi(x)}Y,\qquad v+\fk\mapsto \left.\frac{\mathrm{d}}{\mathrm{d}t}\right|_{t=0}\pi(\exp(tv)x).
$$
Each $j_x$ is a surjective linear map. It is straightforward to see that $$
j_{kx}\circ \Ad_k\big|_{\fg/\fk}=j_x\;(k\in K);
\qquad
j_{f(x)}\circ A=DF_{\pi(x)}\circ j_x.
$$
For each $h\in\cZ^{\Lip}(F)$, we define a map $\cA^{(h)}\colon X\to W$ as follows: for any $x\in X$, let $\cA^{(h)}_x\in W$ denote the linear map $$
\cA^{(h)}_x\colon E_0^*\to \Sym^2((\fg/\fk)^*),\qquad R\mapsto j_x^*(h^*R)_{\pi(x)}.
$$
Since $h$ is Lipschitz, we see that the correspondence $x\mapsto \cA^{(h)}_x$ is bounded and measurable. It is straightforward to see that each $\cA^{(h)}$ satisfies \eqref{E:a-t-1} and \eqref{E:a-t-2}. Then it follows from Lemma \ref{L:regularity} that the set $$
\left\{\cA^{(h)}:h\in\cZ^{\Lip}(F)\right\}
$$ 
is contained in a finite-dimensional space where each member has a smooth representative. We write this property as (FS).

In particular, since $E_0^*$ is finite-dimensional, we see that the collection $$
\left\{
x\mapsto j_x^*(h^*R)_{\pi(x)}:h\in\cZ^{\Lip}(F),\;R\in E_0^*
\right\}
$$
also has property (FS). Moreover, since the correspondence $x\mapsto j_x^*$ is smooth, and each linear map $j_x^*$ is injective, the following vector space also has property (FS): \begin{equation}
\cH:=\Span_{\RR}\left\{h^*R:h\in\cZ^{\Lip}(F),\;R\in E_0^*\right\}.
\end{equation}
\smallskip

\textsl{Step 2: An embedding of the Lipschitz centralizer into a group of isometries.}

Consider the following nonempty open subset of $\cH$: $$
\cC:=\left\{
R\in \cH:R\text{ is positive-definite}
\right\}.
$$
We define a Riemannian metric on $\cC$: for any $R\in \cC$ and any $U,V\in T_R\cC\cong \cH$, put $$
\gamma_R(U,V):=\int_Y\tr(R^{-1}UR^{-1}V)\mathrm{d}\vol_{R},
$$
This gives a Riemannian metric on the product space $Z:=Y\times \cC$: for any $(y,R)\in Z$ and any $(u,U),(v,V)\in T_yY\times T_R\cC$, put \begin{equation}\label{E:prod-metric}
G_{(y,R)}((u,U),(v,V)):=R_y(u,v)+\gamma_R(U,V).
\end{equation}

For any $h\in\cZ^{\Lip}(F)$, it is straightforward to see that  \begin{equation*}
((h^{-1})^*R)_{h(y)}(Dh(u),Dh(v))=R_y(u,v),\qquad
\gamma_{h^*R}(h^*U,h^*V)=\gamma_R(U,V).
\end{equation*}
Write $$
\widehat{h}\colon Z\to Z,\qquad (y,R)\mapsto (h(y),(h^{-1})^*R).
$$
Then each $\widehat{h}$ is an isometry on $Z$ with respect to the Riemannian metric \eqref{E:prod-metric}. Moreover, there is a continuous injective homomorphism: \begin{equation}\label{E:z-i-embed}
\cZ^{\Lip}(F)\to \Isom(Z),\qquad h\mapsto \widehat{h}.
\end{equation}
It is clear that the map \eqref{E:z-i-embed} has a closed image. By \cite[Theorem 8 and Theorem 10]{MS39}, we conclude that $
\cZ^{\Lip}(F)=\cZ^{\infty}(F)$ is a finite-dimensional Lie group. The dimension bound follows from \cite[Proposition 6]{DWWX}.
\end{proof}

Now we proceed to the proof of the ergodic case of Theorem \ref{T:bi-homo}.

\begin{proof}[Proof of the ergodic case of Theorem \ref{T:bi-homo}] In view of Proposition \ref{P:Lip-smooth}, it suffices to find a finite-dimensional real vector space $E_0^*\subseteq\Gamma^{\infty}(\Sym^2T^*Y)$ such that $F^*(E_0^*)=E_0^*$ and $E_0^*$ contains a Riemannian metric. 

Write $A=\Ad_a\circ D\Phi\in\Aut(\fg)$. For $x\in X$, define $$
j_x\colon \fg\to T_{\pi(x)}Y,\qquad v\mapsto \left.\frac{\mathrm{d}}{\mathrm{d}t}\right|_{t=0}\pi(\exp(tv)x).
$$
It is straightforward to see that \begin{equation}\label{E:j-inv}
j_{kx}\circ \Ad_k=j_x\;(k\in K);
\qquad
j_{f(x)}\circ A=DF_{\pi(x)}\circ j_x.
\end{equation}

For any $q\in (\Sym^2\fg)^{\Ad(K)}$, we extend it to a global smooth symmetric contravariant $2$-tensor as $$
Q^{(q)}\colon X\to \Sym^2TY,\qquad x\mapsto (j_x\otimes j_x)q.
$$
It follows from \eqref{E:j-inv} that each $Q^{(q)}$ is left $K$-invariant and hence descends to an element in $\Gamma^{\infty}(\Sym^2TY)$. Note that $F_*Q^{(q)}=Q^{((A\otimes A)q)}$. Then $$
E_0:=\{Q^{(q)}:q\in (\Sym^2\fg)^{\Ad(K)}\}
$$
is an $F_*$-invariant finite-dimensional subspace of $\Gamma^{\infty}(\Sym^2TY)$. In particular, when $q$ is positive-definite, we see that $Q^{(q)}$ is a cometric. 

For any $Q\in \Gamma^{\infty}(\Sym^2TY)$, let $\det_m(Q)$ denote the determinant of $Q$ relative to the fixed volume $m$, and define the smooth symmetric covariant $2$-tensor $\adj_m(Q)$ by the density-normalized cofactor formula. When $Q$ is invertible, this gives
\[
\adj_m(Q)=(\det_m Q)\cdot Q^{-1}.
\]
We note that each $\adj_m(Q)$ is a polynomial of degree $(d-1)$ in $Q$, where $d=\dim{Y}$. Moreover, we have $F^*\adj_m(Q)=\adj_m(F^{-1}_*Q)$. Then $$
E_0^*:=\Span\{\adj_m(Q):Q\in E_0\}
$$
is an $F^*$-invariant finite-dimensional subspace of $\Gamma^{\infty}(\Sym^2T^*Y)$. In particular, when $Q$ is a cometric, we see that $\adj_m(Q)$ is a Riemannian metric. 
\end{proof}

\section{Application (II): Frame flows over negatively curved manifolds}\label{S:fram-flows}
In this section, we prove Theorem \ref{T:frame-flow} as another analogue of Theorem \ref{T:main} in the setting of time-$T$ maps of frame flows over negatively curved manifolds.

For a connected, closed, oriented Riemannian manifold $N$ with $n=\dim{N}$, its full frame bundle $FN$ is a $\SO(n-1)$-principal bundle over its unit tangent bundle $SN$. Write $\pi$ to be the projection map and let $K=\SO(n-1)$. The frame flow $\Phi=\{\Phi_t:t\in\RR\}$ on $FN$ is an isometric extension of the geodesic flow on $SN$.

\subsection{The non-ergodic case} To prove the non-ergodic case of Theorem \ref{T:frame-flow}, we start with the existence of a nonconstant $\Phi_T$-invariant smooth function on $FN$ whose restriction to each $\pi$-fiber is nontrivial. 

\begin{lemma}\label{L:frame-flow-inv-function}
    If $\Phi_T$ is non-ergodic, then there exists a nonconstant $\Phi$-invariant smooth function $f:FN\to [0,1]$ whose restriction to every fiber of $\pi$ is surjective. 
\end{lemma}
\begin{proof}
    Since an ergodic frame flow is mixing (see \cite{Brin75, Lef23}), we see that the non-ergodicity of $\Phi_T$ implies the non-ergodicity of $\Phi=\{\Phi_t:t\in\RR\}$. 
    Then it follows from \cite[Corollary 3.11]{Lef23} that there is a proper closed subgroup $H\leq K$ and a smooth $\Phi$-invariant principal $H$-subbundle $Q\subseteq FN$ over $SN$, on which the restricted flow is ergodic for its natural measure. We define $$
    F\colon FN\to H\backslash K,\qquad F(qk)=Hk\quad (q\in Q,\;k\in K).
    $$
    It is straightforward to see that $F$ is a well-defined nonconstant $\Phi$-invariant smooth map. The conclusion follows by composing $F$ with a surjective smooth function $H\backslash K\to [0,1]$.
\end{proof}

Now we proceed to the proof of the non-ergodic case of Theorem \ref{T:frame-flow}.

\begin{proof}[Proof of the non-ergodic case of Theorem \ref{T:frame-flow}] By Lemma \ref{L:frame-flow-inv-function}, let $f:FN\to [0,1]$ be a nonconstant $\Phi$-invariant smooth function whose restriction to every fiber of $\pi$ is surjective. For each $\rho\in C_c^{\infty}((0,1),\RR)$, write $$
h^{(\rho)}\colon FN\to FN,\qquad u\mapsto \Phi_{\rho(f(u))}(u).
$$
It is clear that $h^{(\rho)}\in \cZ^{\infty}(\Phi)$. Moreover, this gives a continuous group homomorphism: $$
C_c^{\infty}((0,1),\RR)\to \cZ^{\infty}(\Phi),\qquad \rho\mapsto h^{(\rho)}. 
$$
We only need to check that it is injective. In fact, suppose that $h^{(\rho)}=\id_{FN}$. Then projecting it to $SN$ gives $\varphi_{\rho(f(u))}(\pi(u))=\pi(u)$ for all $u\in FN$. Since $(N,g)$ has negative curvature, there exists a non-periodic geodesic on $N$, say $\gamma_{\xi}$ where $\xi\in SN$. It follows that $\rho(f(u))=0$ for any $u\in \pi^{-1}(\xi)$. Since $f$ is surjective onto $[0,1]$ on any $\pi$-fiber, we conclude that $\rho\equiv 0$. This completes the proof.
\end{proof}

\subsection{The ergodic case} Now we prove the ergodic case using the expansiveness of the geodesic flow on $SN$.

\begin{proof}[Proof of the ergodic case of Theorem \ref{T:frame-flow}] 
Since $\Phi_t$ commutes with the right action of $K$, the following map is a continuous group homomorphism: \begin{equation}\label{E:centralizer-r-so}
    \RR\times K\to \cZ^0(\Phi_T),\qquad (t,k)\mapsto \Phi_t(\cdot)k^{-1}.
\end{equation}

We first show that the map \eqref{E:centralizer-r-so} is injective. Suppose that $\Phi_t(\cdot)k^{-1}=\Phi_{t'}(\cdot)(k')^{-1}$. Applying $\pi$ to both sides of this identity gives $\varphi_t\circ \pi=\varphi_{t'}\circ \pi$. Since a negatively curved manifold has no conjugate points, we have $t=t'$. Then the freeness of $K$-action gives $k=k'$. This verifies the injectivity.

We next show that the image of \eqref{E:centralizer-r-so} contains a small neighborhood of $\id_{FN}$ in $\cZ^0(\Phi_T)$. In fact, since $SN$ is compact, for some sufficiently small $\epsilon>0$ the map \begin{equation}\label{E:compact-homeo}
    SN\times [-\epsilon,\epsilon]\to SN\times SN,\qquad (\xi,t)\mapsto (\xi,\varphi_t(\xi))
\end{equation}
is a homeomorphism onto its image. Since the geodesic flow $\{\varphi_t:t\in\RR\}$ on $SN$ is expansive (see \cite[Corollary 5.3.5]{FH19}), there exists some $\delta=\delta(\epsilon)>0$ such that for any $\xi,\xi'\in SN$, \begin{equation}\label{E:expansive}
\sup_{t\in\RR}\dist_{SN}(\varphi_t(\xi),\varphi_t(\xi'))<\delta\;\Longrightarrow\; \xi=\varphi_{\tau}(\xi')\text{ for some unique }\tau\in\RR\text{ with }|\tau|<\epsilon.
\end{equation}
Write $$
C_{\pi}\coloneqq\sup_{u\in FN}\|D_u\pi\|,\qquad C_{T}\coloneqq\sup_{\substack{u\in FN\\ |r|\leq |T|/2}}\|D_u\Phi_r\|.
$$
We claim that: the image of \eqref{E:centralizer-r-so} contains the following neighborhood of $\id_{FN}$:
$$\cU:=\left\{h\in \cZ^0(\Phi_T):\dist_{C^0}(h,\id_{FN})<C_{\pi}^{-1}C_T^{-1}\delta\right\}.$$

In fact, let $h\in \cU$ and $u\in FN$. For any $t\in\RR$, we choose $j\in\ZZ$ and $|r|\leq |T|/2$ with $t=jT+r$. Then we have \begin{equation*}
    \begin{aligned}
        \dist_{SN}(\varphi_t\circ \pi(h(u)),\varphi_t\circ \pi(u))
        &=\dist_{SN}(\pi\circ\Phi_t(h(u)),\pi\circ\Phi_t(u))\\
        &=\dist_{SN}(\pi\circ\Phi_r\circ h(\Phi^j_T(u)),\pi\circ\Phi_r\circ\Phi^j_T(u))\\
        &\leq C_{\pi}C_{T}\cdot \dist_{FN}(h(\Phi^j_T(u)),\Phi^j_T(u))\\
        &\leq C_{\pi}C_{T}\cdot \dist_{C^0}(h,\id_{FN})<\delta.
    \end{aligned}
\end{equation*}
It follows from \eqref{E:expansive} that there exists a unique $\tau_h(u)\in\RR$ with $|\tau_h(u)|<\epsilon$ such that \begin{equation}\label{E:def-tau}
\pi(h(u))=\varphi_{\tau_h(u)}(\pi(u)).
\end{equation}
By \eqref{E:compact-homeo}, we see that $\tau_h(u)$ depends continuously in $h\in \cU$ and $u\in FN$. Since the right action of $K$ on each $\pi$-fiber is transitive and free, we see from \eqref{E:def-tau} that there exists a unique $\kappa_h(u)\in K$ with \begin{equation}\label{E:def-kappa}
    h(u)=\Phi_{\tau_h(u)}(u)\cdot \kappa_h(u)^{-1}. 
\end{equation}
Moreover, $\kappa_h(u)$ also depends continuously in $h\in \cU$ and $u\in FN$.
Then we have \begin{equation*}
\Phi_{T+\tau_h(u)}(u)\cdot \kappa_h(u)^{-1}=\Phi_T\circ h(u)=h\circ \Phi_T(u)=\Phi_{\tau_h(\Phi_T(u))+T}(u)\cdot \kappa_h(\Phi_T(u))^{-1}. 
\end{equation*}
Applying $\varphi_{-T}\circ\pi$ on both sides gives $\varphi_{\tau_h(u)}(\pi(u))=\varphi_{\tau_h(\Phi_T(u))}(\pi(u))$. In view of \eqref{E:compact-homeo}, we see that $\tau_h(u)=\tau_h(\Phi_T(u))$, and hence that $\kappa_h(u)=\kappa_h(\Phi_T(u))$. Since $\Phi_T$ is ergodic, we conclude that both $\tau_h(\cdot)$ and $\kappa_h(\cdot)$ are constant on $FN$. Therefore, we see from \eqref{E:def-kappa} that there exists some $\tau\in\RR$ and $k\in K$ such that  $$
h(u)=\Phi_{\tau}(u)k^{-1}\quad \text{ for all }u\in FN.
$$
This verifies the claim.

Now we conclude that the map \eqref{E:centralizer-r-so} is an injective continuous group homomorphism whose image contains a small neighborhood of $\id_{FN}$ in $\cZ^0(\Phi_T)$. So by connectedness, its image equals the identity component of $\cZ^0(\Phi_T)$. Since the image of \eqref{E:centralizer-r-so} consists of smooth diffeomorphisms, we see that $\cZ^0(\Phi_T)$ and $\cZ^{\infty}(\Phi_T)$ have the same identity component, which is isomorphic to $\RR\times K$. The proof is complete.
\end{proof}

\section{Application (III): Automorphisms of Kummer-type complex surfaces}\label{S:Kummer}
In this section, we prove Theorem \ref{T:c-surface} as another analogue of Theorem \ref{T:main} in the setting of automorphisms of Kummer-type complex surfaces.

For later use, we fix some notation. For a compact torus $T=\CC^2/\Gamma$, let $G$ be a nontrivial finite subgroup of $\Aut(\CC^2,\Gamma)$, and let $f\in\Aut(\CC^2,\Gamma)$ normalize $G$. Consider the quotient map and the resolution map $$
q\colon T\to Y=T/G,\quad\hbox{ and }\quad \pi\colon X\to Y.
$$
Let $Y_{\reg}$ denote the complement of finitely many singularities in $Y$, and write $T_{\reg}=q^{-1}(Y_{\reg})$, $X_{\reg}=\pi^{-1}(Y_{\reg})$. Then $q\big|_{T_{\reg}}$ is a $|G|$-sheeted covering, and $\pi\big|_{X_{\reg}}$ is a bi-holomorphism. It follows that the map \begin{equation}\label{E:reg-cover}
p:=(\pi^{-1}\circ q)\big|_{T_{\reg}}\colon T_{\reg}\to X_{\reg}
\end{equation}
is a holomorphic $|G|$-sheeted covering sending a Haar ($m_T$)-full subset in $T$ to a $\mu_X$-full subset in $X$. Thus, the system $(X,\mu_X,F)$ is a measurable factor of $(T,m_T,f)$.

We first prove the equivalence ``(1) $\Longleftrightarrow$ (2)'' in Theorem \ref{T:c-surface}. This is a simple consequence of the construction in Lemma \ref{L:torus}.

\begin{lemma}\label{L:hyper-erg}
    $f$ is hyperbolic if and only if $(X,\mu_X,F)$ is ergodic (or mixing).
\end{lemma}
\begin{proof} First, suppose that $f$ is hyperbolic. Since a hyperbolic toral automorphism is mixing, we see that its measurable factor $(X,\mu_X,F)$ is also mixing, and in particular, it is ergodic.

Conversely, suppose that $f$ is not hyperbolic. Since $f\in\Aut(\CC^2,\Gamma)$, it must be quasi-unipotent. We write $f^m=I+N$ for some $m\geq 1$ and some $N$ with $N^2=0$.  It follows that the induced dynamics on the nontrivial factor $T/\langle \mathrm{Im}(N),G\rangle$ is periodic. So $(X,\mu_X,F)$ is non-ergodic.
\end{proof}

\subsection{The non-ergodic case} 
In this subsection, we prove the implication ``(3)$\Longrightarrow$(2)'' of Theorem \ref{T:c-surface}, which follows from the statement in the non-ergodic case.

\begin{proof}[Proof of the non-ergodic case of Theorem \ref{T:c-surface}]
In view of Proposition \ref{P:function-vector-field}, it suffices to show the existence of the followings: \begin{enumerate}
    \item[(i)] a nonconstant $F$-invariant smooth function $\Psi$ on $X$;
    \item[(ii)] an $F$-invariant invariant smooth vector field $\fX$ on $X$, and a small open neighborhood $\cU\subseteq X$ in which every point $x$ satisfies $$
\fX_x\neq 0\quad\text{and}\quad (\mathrm{d}\Psi)_x\neq 0.
$$
\end{enumerate}

We keep the notations as in the proof of Lemma \ref{L:hyper-erg}. Since $(X,\mu_X,F)$ is not ergodic, it follows that the induced dynamics of $f$ on the nontrivial factor $\overline{T}:=T/\mathrm{Im}(N)$ is periodic. We choose any nonzero, nonnegative smooth function $\rho$ on $\overline{T}$ vanishing on a small open neighborhood of the image of $T\setminus T_{\reg}$. By averaging over the $f$-action and the $G$-action, we may assume that $\rho$ is $f$-invariant and $G$-invariant. Then it descends to a nonconstant $F$-invariant smooth function on $T/\langle \mathrm{Im}(N),G\rangle$. Lifting it back to $X$ gives a nonconstant $F$-invariant smooth function $\Psi$ on $X$ vanishing on a small open neighborhood of $X\setminus X_{\reg}$. This completes (i).

We also choose an $f$-invariant and $G$-invariant Riemannian metric on $\overline{T}$. Then the gradient $\nabla\rho$ is a nonzero $f$-invariant, $G$-invariant, smooth vector field on $\overline{T}$. Since $N^2=0$, we define the linear isomorphism $S\colon T_0\overline{T}\to \ker(DN)\subseteq T_0T$ by $$
S:=\begin{cases}
    \id:T_0T\to T_0T, &\text{if }N=0;\\
    DN:T_0\overline{T}=T_0T/\mathrm{Im}(DN)\to \mathrm{Im}(DN), &\text{if }N\neq 0.
\end{cases}
$$
Note that in the second case we have $\ker(DN)=\mathrm{Im}(DN)$. We identify all tangent spaces of $T$ and $\overline{T}$ with their tangent spaces at the origin by translations, and define
$$\fT_z:=S\bigl(\nabla\rho_{z+\mathrm{Im}(N)}\bigr),\qquad \text{for all }z\in T.$$
It follows that $\fT$ is a nonzero $f$-invariant, $G$-invariant, smooth vector field on $T$ vanishing on a neighborhood of $T\setminus T_{\reg}$. Under the covering map $p$, it descends to a nonzero $F$-invariant smooth vector field $\fX$ on $X$ vanishing on a neighborhood of $X\setminus X_{\reg}$. When $\mathrm{d}\rho_{z+\mathrm{Im}(N)}\neq 0$, we have $(\mathrm{d}\Psi)_{p(z)}\neq 0$ and $\fX_{p(z)}\neq 0$. This completes (ii).
\end{proof}

\subsection{The ergodic case} 
Now we prove the implication ``(2)$\Longrightarrow$(4)'' in Theorem \ref{T:c-surface}. This is based on some arguments in \cite[Proof of Theorem B]{DWWX} and \cite[Lemmas 3.1 and 3.2]{Zhang26}, as well as a local boundedness argument near the exceptional divisor.

\begin{lemma}
    If $(X,\mu_X,F)$ is ergodic, then only essentially bounded $F$-invariant almost complex structures on $X$ are $\pm J_0$. 
\end{lemma}
\begin{proof}
We divide the proof into several steps.\medskip

\textsl{Step 1: Lift an invariant almost complex structure to a constant endomorphism.}

Let $J$ be an essentially bounded $F$-invariant almost complex structure on $X$. This means that $J$ is an essentially bounded section of the real bundle $\End_{\RR}(TX)\to X$, satisfying that $$
J_x^2=-\id_{T_xX},\quad\text{and}\quad DF_x\circ J_x=J_{F(x)}\circ DF_x,\quad\text{for }\mu_X\text{-almost all }x\in X. 
$$
For any $z\in T_{\reg}$, write $$
B_z:=Dp_z^{-1}\circ J_{p(z)}\circ Dp_z\in \End_{\RR}(T_zT).
$$
It is straightforward to see that $B_{f(z)}=Df_z\circ B_z\circ Df_z^{-1}$. By identifying each tangent space $T_zT$ with $\RR^4$ and defining $B_z$ arbitrarily for $z\in T\setminus T_{\reg}$, we regard $z\mapsto B_z$ as a map $T\to \End_{\RR}(\RR^4)$. Then we claim that: \begin{equation}\label{E:B-ess-bounded}
    B\colon T\to \End_{\RR}(\RR^4),\quad z\mapsto B_z\quad\text{is essentially bounded}.
\end{equation}

In fact, under the above identification $T_zT\cong\RR^4$, we consider a linear map $$
L:\End_{\RR}(\RR^4)\to \End_{\RR}(\RR^4),\qquad C\mapsto Df\circ C\circ Df^{-1}.
$$
It is straightforward to see that $B_{f(z)}=L(B_z)$ for $m_T$-almost all $z\in T$. Write $$
W:=\bigoplus_{|\lambda|=1}\ker(L-\lambda\cdot \id).
$$
We first verify that: for $m_T$-almost all $z\in T$, one has $B_z\in W$. Indeed, we fix any norm $\|\cdot\|$ on $\End_{\RR}(\RR^4)$, and for any $R\geq 1$, write $E_R:=\{z\in T_{\reg}:\|B_z\|\leq R\}$. It is clear that $m_T\left(\bigcup_{R\geq 1}E_R\right)=1$. For each $R\geq 1$ and $m_T$-almost every $z\in E_R$, by Poincare recurrence, there exist two sequences $n_j\to +\infty$ and $n_j'\to +\infty$ such that $$
f^{n_j}(z)\in E_R,\quad\text{and}\quad f^{-n'_j}(z)\in E_R,\quad \text{for all }j\geq 1.
$$
that is, $\|L^{n_j}(B_z)\|\leq R$ and $\|L^{-n'_j}(B_z)\|\leq R$ for all $j\geq 1$. This means that each $B_z$ must lie in the direct sum of all generalized unit-modulus eigenspaces of $L$. Since $f$ is hyperbolic, we see that $Df$ is $\CC$-diagonalizable, so is $L$. Thus, we conclude that $B_z\in W$ for $m_T$-almost all $z\in T$.

Now we choose and fix a norm on $W$ that is $L$-invariant. Then we see that for $m_T$-almost all $z\in T$, one has $\|B_{f(z)}\|=\|L(B_z)\|=\|B_z\|$. By the ergodicity of $f$ on $T$, the map $z\mapsto \|B_z\|$ is $m_T$-almost surely constant. This proves our claim \eqref{E:B-ess-bounded}.

Finally, since $f$ is mixing on $T$ (see Lemma \ref{L:hyper-erg}), by \cite[Lemma 3.2]{Zhang26}, we see that there exists a constant $B_0\in \End_{\RR}(\RR^4)$ with $B_z=B_0$ for $m_T$-almost all $z\in T$ and $B_0=Df\circ B_0\circ Df^{-1}$. 

Moreover, for $\gamma\in G$ and $z\in T_{\reg}$, since $p(\gamma z)=p(z)$, we see that $B_{\gamma z}=\gamma\circ B_z\circ \gamma^{-1}$, and hence that $B_0=\gamma\circ B_0\circ \gamma^{-1}$. We also see from $J_{p(z)}^{2}=-\id$ that $B_0^2=-\id$.
\medskip

\textsl{Step 2: Local boundedness of the complex structure near the singularity.}

Since $f$ is hyperbolic on a complex torus of dimension two, we can choose complex coordinates $(z_1,z_2)$ along its stable and unstable eigenspaces such that
$$
Df=\begin{pmatrix}
\alpha&0\\
0&\beta
\end{pmatrix},
\qquad\text{where } |\alpha|<1<|\beta|.
$$
Moreover, we can also diagonalize $G$ in these coordinates. Indeed, since $f$ normalizes the finite group $G$, some power $f^m$ centralizes $G$. Then every element of $G$ preserves the stable and unstable eigenspaces of $Df^m$, which are the same as those of $Df$.

Since $G\leq\SL_2(\mathbb C)$, it follows that $G$ is cyclic. We write
$$G=\langle\gamma\rangle,\qquad\gamma(z_1,z_2)=(\zeta z_1,\zeta^{-1}z_2),$$
where \(\zeta\) is a primitive \(r\)-th root of unity and
$r=|G|\ge2.$ Then a singular point of $Y=T/G$ has its local model as
$\mathbb C^2/\langle\gamma\rangle.$ A chart of the minimal resolution of $\mathbb C^2/\langle\gamma\rangle$ is given by
$$w_1=z_1^r,\qquad w_2=\frac{z_2}{z_1^{r-1}}.$$
We calculate that (in a local branch for $z_1=w_1^{1/r}$) \begin{equation}\label{E:D-z-w}
D_{(z_1,z_2)}(w_1,w_2)=\begin{pmatrix}
rz_1^{r-1}&0\\
-(r-1)z_2z_1^{-r}&z_1^{-(r-1)}
\end{pmatrix}=
\begin{pmatrix}
rw_1^{(r-1)/r}&0\\
-(r-1)w_2w_1^{-1/r}&w_1^{-(r-1)/r}
\end{pmatrix}.
\end{equation}

With respect to the complex coordinates $z=(z_1,z_2)\in\CC^2$, we write the map $B_0\in\End_{\RR}(\RR^4)=\End_{\RR}(\CC^2)$ as
$$B_0(z)=Pz+Q\overline{z}\qquad\text{for all }z\in\CC^2,\qquad\text{where }P,Q\in \CC^{2\times 2}.
$$
Then $B_0\circ Df=Df\circ B_0$ is equivalent to
$$
P\circ Df=Df\circ P,
\quad\text{and}\quad 
Q\circ \overline{Df}=Df\circ Q.
$$
Since $|\alpha|\neq |\beta|$, these force both \(P\) and \(Q\) to be diagonal:
$$
P=\diag(a_1,a_2),
\qquad
Q=\diag(b_1,b_2).
$$
Moreover, with respect to the resolution coordinates $w=(w_1,w_2)\in\CC^2$, we write the induced map $p_*B_0$ as $$
(p_*B_0)(w)=\widetilde{P}w+\widetilde{Q}\overline{w}\qquad\text{for all }w\in\CC^2,
$$
where for $D:=D_{(z_1,z_2)}(w_1,w_2)$ and $\rho:=(w_1/\overline{w_1})^{1/r}$, \begin{subequations}
    \begin{align}
        &{
        \widetilde{P}=DPD^{-1}=\begin{pmatrix}
    a_1 &0 \\ \dfrac{(r-1)(a_2-a_1)w_2}{rw_1} &a_2
\end{pmatrix},}\label{E:tilde-P}\\
&{
\widetilde{Q}=DQ\overline{D^{-1}}=
\begin{pmatrix}
b_1\rho^{\,r-1} &0\\
\dfrac{r-1}{r}\left(\dfrac{b_2\overline{w_2}\rho}{w_1}-\dfrac{b_1w_2}{\rho\overline{w_1}}\right)
&b_2\rho^{-(r-1)}
\end{pmatrix}.
}\label{E:tilde-Q}
    \end{align}
\end{subequations}
We claim that: $P$ is a scalar and $Q=0$. 

In fact, since $p_*B_0=J$ is essentially bounded, we see that $a_1=a_2$ and $b_1=b_2=0$, otherwise the $(2,1)$-entry of $\widetilde{P}$ and $\widetilde{Q}$ would be unbounded as $w_1\to 0$. This means that $P$ is a scalar and $Q=0$. 

Therefore, we conclude that $B_0=a_1\cdot\id_{\CC^2}$. Since $B_0^2=-\id$, we have $a_1=\pm i$, that is, $B_0$ is one of the two standard complex structures on $T$. It follows that $J$ is one of the two standard complex structures on $X$.
\end{proof}

Finally, we prove the implication ``(4)$\Longrightarrow$(3)'' in Theorem \ref{T:c-surface}.

\begin{lemma}
Let $F$ be a holomorphic automorphism of a K3 surface $(X,J_0)$. Suppose that any essentially bounded invariant almost complex structure equals $\pm J_0$ almost everywhere. Then $\cZ^{\Lip}(F)=\cZ^{\infty}(F)\subseteq\Aut^{\pm}(X)$ is a discrete group.
\end{lemma}
\begin{proof}
    Let $h\in\cZ^{\Lip}(F)$. Then $h^*J_0$ is an essentially bounded invariant almost complex structure. By our assumption, we see that $h^*J_0=\pm J_0$, that is, $Dh$ is complex linear or complex anti-linear. In view of \cite[Lemma 3.3]{Zhang26}, this implies that $h$ is holomorphic or anti-holomorphic. In particular, we see that $\cZ^{\Lip}(F)=\cZ^{\infty}(F)$ is a subgroup of $\Aut^{\pm}(X)$. Finally, the discreteness of $\Aut^{\pm}(X)$ follows from \cite[Chapter~15]{Huybrechts}.
\end{proof}

\end{document}